\documentclass[12pt]{article}
\PassOptionsToPackage{hyphens}{url}
\usepackage[pagebackref=true,colorlinks]{hyperref}
\usepackage{amsmath,amssymb,amsthm,mathtools,amsrefs,mathrsfs}
\usepackage{pifont}
\usepackage{bbm}
\usepackage{tikz}
\usepackage{adjustbox}
\usepackage{stmaryrd}
\usepackage{fancyhdr}
\usepackage{dsfont}
\usepackage[normalem]{ulem}
\usepackage[bottom]{footmisc}
\usepackage{enumerate}
\usepackage[all]{xy} 
\usetikzlibrary{circuits.ee.IEC}
\hypersetup{
    colorlinks,
    linkcolor={red!50!black},
    citecolor={blue!50!black},
    urlcolor={blue!80!black}
}

\usepackage{iftex}

\ifpdf
  \usepackage[T1]{fontenc}        
\else
  \usepackage{breakurl}           
\fi

\usepackage{tocbasic}
\DeclareTOCStyleEntry[
  beforeskip=0.3em plus 1pt,
  pagenumberformat=\textbf
]{tocline}{section}

\makeatletter
\newcommand*{\shifttext}[2]{%
  \settowidth{\@tempdima}{#2}%
  \makebox[\@tempdima]{\hspace*{#1}#2}%
}
\makeatother
\makeatletter
\renewcommand*\env@matrix[1][\arraystretch]{%
  \edef\arraystretch{#1}%
  \hskip -\arraycolsep
  \let\@ifnextchar\new@ifnextchar
  \array{*\c@MaxMatrixCols c}}
\makeatother

\theoremstyle{plain}
\newtheorem{theorem}[equation]{Theorem}
\newtheorem{lemma}[equation]{Lemma}
\newtheorem{proposition}[equation]{Proposition}
\newtheorem{corollary}[equation]{Corollary}
\theoremstyle{definition}
\newtheorem{definition}[equation]{Definition}
\newtheorem{construction}[equation]{Construction}
\newtheorem{question}[equation]{Question}

\newtheorem{problem}[equation]{Problem}
\newtheorem{example}[equation]{Example}
\newtheorem{exercise}[equation]{Exercise}
\newtheorem*{answer}{Answer}
\newtheorem*{solution}{Solution}
\newtheorem{remark}[equation]{Remark}

\newtheorem{notation}[equation]{Notation}
\newtheorem{noterm}[equation]{Notation and Terminology}
\newtheorem{term}[equation]{Terminology}

\newcommand\define[1]{\emph{\textbf{#1}}}

\numberwithin{equation}{section}

\let\C=\Chi

\newcommand{\be}{\begin{equation}}
\newcommand{\ee}{\end{equation}}
\def\ba{\begin{align}} 
\def\ea{\end{align}}
\newcommand{\bea}{\begin{eqnarray}}
\newcommand{\eea}{\end{eqnarray}}
\newcommand{\bx}{\begin{example}}
\newcommand{\ex}{\end{example}}
\newcommand{\bex}{\begin{exercise}}
\newcommand{\eex}{\end{exercise}}
\newcommand{\ban}{\begin{answer}}
\newcommand{\ean}{\end{answer}}
\newcommand{\bt}{\begin{theorem}}
\newcommand{\et}{\end{theorem}}
\newcommand{\bc}{\begin{corollary}}
\newcommand{\ec}{\end{corollary}}
\newcommand{\blem}{\begin{lemma}}
\newcommand{\elem}{\end{lemma}}
\newcommand{\bp}{\begin{problem}}
\newcommand{\ep}{\end{problem}}
\newcommand{\bn}{\begin{proposition}}
\newcommand{\en}{\end{proposition}}
\newcommand{\bd}{\begin{definition}}
\newcommand{\ed}{\end{definition}}
\newcommand{\bcon}{\begin{construction}}
\newcommand{\econ}{\end{construction}}
\newcommand{\bq}{\begin{question}}
\newcommand{\eq}{\end{question}}
\newcommand{\bprf}{\begin{proof}}
\newcommand{\eprf}{\end{proof}}
\newcommand{\br}{\begin{remark}}
\newcommand{\er}{\end{remark}}
\newcommand{\bs}{\begin{solution}}
\newcommand{\es}{\end{solution}}
\newcommand{\beqs}{\begin{eqnarray}}
\newcommand{\eeqs}{\end{eqnarray}}
\newcommand{\bnt}{\begin{noterm}}
\newcommand{\ent}{\end{noterm}}
\newcommand{\bnot}{\begin{notation}}
\newcommand{\enot}{\end{notation}}

 \let\ov=\overline

\newcommand{\<}{\langle}
\renewcommand{\>}{\rangle}

\newcommand{\id}{\mathrm{id}}

\newcommand{\Tr}{\operatorname{Tr}}

\newcommand{\map}[1]{\mathcal{#1}}

\def\C{{{\mathbb C}}}

\newcommand\aeequals[1]{\underset{\raisebox{0.3ex}[0pt][0pt]{\scriptsize${#1}$}}{=}}

\newcommand{\Ad}{\mathrm{Ad}}
\def\Hi{{{\mathcal{H}}}}

\DeclareFontFamily{OT1}{pzc}{}
\DeclareFontShape{OT1}{pzc}{m}{it}{ <-> s*[1.2] pzcmi7t }{}
\DeclareMathAlphabet{\mathpzc}{OT1}{pzc}{m}{it}
\newcommand{\Alg}[1]{\mathpzc{#1}}

\newcommand{\matr}{\mathbb{M}}

\newcommand{\op}{\mathrm{op}}

\newcommand{\FinStoch}{\mathbf{FinStoch}}
\newcommand{\NCFinStoch}{\mathbf{NCFinStoch}}
\newcommand{\FinProb}{\mathbf{FinProb}}

\def\retro{\mathscr{R}}
\def\Petz{\mathrm{P}}

\def\Bayes{\mathrm{B}}

\def\states{\mathscr{S}}

\def\stFS{\mathscr{S}^{+}(\FinStoch)}

\newcommand{\ben}{\renewcommand{\theenumi}{\alph{enumi}} 
\renewcommand{\labelenumi}{(\theenumi)}\begin{enumerate}}
\newcommand{\een}{\end{enumerate}}

\newcommand\blfootnote[1]{%
  \begingroup
  \renewcommand\thefootnote{}\footnote{#1}%
  \addtocounter{footnote}{-1}%
  \endgroup
}

\newcommand{\bsm}{\begin{smallmatrix}}
\newcommand{\esm}{\end{smallmatrix}}

\newlength\stateheight
\newlength\minimumstatewidth
\tikzset{width/.initial=\minimummorphismwidth}
\tikzset{colour/.initial=white}

\newif\ifblack\pgfkeys{/tikz/black/.is if=black}
\newif\ifwedge\pgfkeys{/tikz/wedge/.is if=wedge}
\newif\ifvflip\pgfkeys{/tikz/vflip/.is if=vflip}
\newif\ifhflip\pgfkeys{/tikz/hflip/.is if=hflip}
\newif\ifhvflip\pgfkeys{/tikz/hvflip/.is if=hvflip}

\pgfdeclarelayer{foreground}
\pgfdeclarelayer{background}
\pgfsetlayers{background,main,foreground}

\def\thickness{0.4pt}

\makeatletter
\pgfkeys{%
  /tikz/on layer/.code={
    \pgfonlayer{#1}\begingroup
    \aftergroup\endpgfonlayer
    \aftergroup\endgroup
  },
  /tikz/node on layer/.code={
    \gdef\node@@on@layer{%
      \setbox\tikz@tempbox=\hbox\bgroup\pgfonlayer{#1}\unhbox\tikz@tempbox\endpgfonlayer\pgfsetlinewidth{\thickness}\egroup}
    \aftergroup\node@on@layer
  },
  /tikz/end node on layer/.code={
    \endpgfonlayer\endgroup\endgroup
  }
}
\def\node@on@layer{\aftergroup\node@@on@layer}
\makeatother

\makeatletter
\pgfdeclareshape{state}
{
    \savedanchor\centerpoint
    {
        \pgf@x=0pt
        \pgf@y=0pt
    }
    \anchor{center}{\centerpoint}
    \anchorborder{\centerpoint}
    \saveddimen\overallwidth
    {
        \pgf@x=3\wd\pgfnodeparttextbox
        \ifdim\pgf@x<\minimumstatewidth
            \pgf@x=\minimumstatewidth
        \fi
    }
    \savedanchor{\upperrightcorner}
    {
        \pgf@x=.5\wd\pgfnodeparttextbox
        \pgf@y=.5\ht\pgfnodeparttextbox
        \advance\pgf@y by -.5\dp\pgfnodeparttextbox
    }
    \anchor{A}
    {
        \pgf@x=-\overallwidth
        \divide\pgf@x by 4
        \pgf@y=0pt
    }
    \anchor{B}
    {
        \pgf@x=\overallwidth
        \divide\pgf@x by 4
        \pgf@y=0pt
    }
    \anchor{text}
    {
        \upperrightcorner
        \pgf@x=-\pgf@x
        \ifhflip
            \pgf@y=-\pgf@y
            \advance\pgf@y by 0.4\stateheight
        \else
            \pgf@y=-\pgf@y
            \advance\pgf@y by -0.4\stateheight
        \fi
    }
    \backgroundpath
    {
       \begin{pgfonlayer}{foreground}
        \pgfsetstrokecolor{black}
        \pgfsetlinewidth{\thickness}
        \pgfpathmoveto{\pgfpoint{-0.5*\overallwidth}{0}}
        \pgfpathlineto{\pgfpoint{0.5*\overallwidth}{0}}
        \ifhflip
            \pgfpathlineto{\pgfpoint{0}{\stateheight}}
        \else
            \pgfpathlineto{\pgfpoint{0}{-\stateheight}}
        \fi
        \pgfpathclose
        \ifblack
            \pgfsetfillcolor{black!50}
            \pgfusepath{fill,stroke}
        \else
            \pgfusepath{stroke}
        \fi
       \end{pgfonlayer}
    }
}

\pgfdeclareshape{my ground}
{
  \inheritsavedanchors[from=rectangle ee]
  \inheritanchor[from=rectangle ee]{center}
  \inheritanchor[from=rectangle ee]{north}
  \inheritanchor[from=rectangle ee]{south}
  \inheritanchor[from=rectangle ee]{east}
  \inheritanchor[from=rectangle ee]{west}
  \inheritanchor[from=rectangle ee]{north east}
  \inheritanchor[from=rectangle ee]{north west}
  \inheritanchor[from=rectangle ee]{south east}
  \inheritanchor[from=rectangle ee]{south west}
  \inheritanchor[from=rectangle ee]{input}
  \inheritanchor[from=rectangle ee]{output}
  \anchorborder{\csname pgf@anchor@my ground@input\endcsname}

  \backgroundpath{
    \pgf@process{\pgfpointadd{\southwest}{\pgfpoint{\pgfkeysvalueof{/pgf/outer xsep}}{\pgfkeysvalueof{/pgf/outer ysep}}}}
    \pgf@xa=\pgf@x \pgf@ya=\pgf@y
    \pgf@process{\pgfpointadd{\northeast}{\pgfpointscale{-1}{\pgfpoint{\pgfkeysvalueof{/pgf/outer xsep}}{\pgfkeysvalueof{/pgf/outer ysep}}}}}
    \pgf@xb=\pgf@x \pgf@yb=\pgf@y
    \pgfpathmoveto{\pgfqpoint{\pgf@xa}{\pgf@ya}}
    \pgfpathlineto{\pgfqpoint{\pgf@xa}{\pgf@yb}}
    \pgfmathsetlength\pgf@xa{.5\pgf@xa+.5\pgf@xb}
    \pgfmathsetlength\pgf@yc{.16666\pgf@yb-.16666\pgf@ya}
    \advance\pgf@ya by\pgf@yc
    \advance\pgf@yb by-\pgf@yc
    \pgfpathmoveto{\pgfqpoint{\pgf@xa}{\pgf@ya}}
    \pgfpathlineto{\pgfqpoint{\pgf@xa}{\pgf@yb}}
    \advance\pgf@ya by\pgf@yc
    \advance\pgf@yb by-\pgf@yc
    \pgfpathmoveto{\pgfqpoint{\pgf@xb}{\pgf@ya}}
    \pgfpathlineto{\pgfqpoint{\pgf@xb}{\pgf@yb}}
  }
}
\makeatother

\tikzset{inline text/.style =
  {text height=1.2ex,text depth=0.25ex,yshift=0.5mm}}
\tikzset{arrow box/.style =
  {rectangle,inline text,fill=white,draw,
    minimum height=5mm,yshift=-0.5mm,minimum width=5mm}}
\tikzset{bubble/.style =
  {inner sep=0mm,minimum width=3mm,minimum height=3mm,
    draw,shape=circle,fill=white}}
\tikzset{dot/.style =
  {inner sep=0mm,minimum width=1mm,minimum height=1mm,
    draw,shape=circle}}
\tikzset{white dot/.style = {dot,fill=white,text depth=-0.2mm}}
\tikzset{scalar/.style = {diamond,draw,inner sep=1pt}}
\tikzset{square/.style =
  {inner sep=0mm,minimum width=2mm,minimum height=2mm,
    draw,shape=rectangle}}

\tikzset{star/.style = {dot,fill=white,text depth=-0.2mm}}
\tikzset{copier/.style = {dot,fill,text depth=-0.2mm}}
\tikzset{fakecopier/.style = {square,fill,text depth=-0.2mm}}
\tikzset{discarder/.style = {my ground,draw,inner sep=0pt,
    minimum width=4.2pt,minimum height=11.2pt,anchor=input,rotate=90}}

\tikzset{xshiftu/.style = {shift = {(#1, 0)}}}
\tikzset{yshiftu/.style = {shift = {(0, #1)}}}
\tikzset{scriptstyle/.style={font=\everymath\expandafter{\the\everymath\scriptstyle}}}

\title{Reversing information flow:\\
\Large retrodiction in semicartesian categories}
\author{Arthur J.~Parzygnat
}

\newcommand{\Addresses}{{
 \bigskip
 \footnotesize
 
   A.~J.~Parzygnat, \textsc{%
 The Experimental Study Group, 
Massachusetts Institute of Technology, 
Cambridge, Massachusetts 02139, USA}\par\nopagebreak
 \textit{E-mail address}, A.~J.~Parzygnat: \texttt{arthurjp@mit.edu}
}}

\begin{document}
\emergencystretch 2em

\maketitle
\begin{abstract}  
In statistical inference, retrodiction is the act of inferring potential causes in the past based on knowledge of the effects in the present and the dynamics leading to the present. 
Retrodiction is applicable even when the dynamics is not reversible, and it agrees with the reverse dynamics when it exists, so that retrodiction may be viewed as an extension of inversion, i.e., time-reversal. 
Recently, an axiomatic definition of retrodiction has been made in a way that is applicable to both classical and quantum probability using ideas from category theory. 
Almost simultaneously, a framework for information flow in in terms of semicartesian categories has been proposed in the setting of categorical probability theory. 
Here, we formulate a general definition of retrodiction to add to the information flow axioms in semicartesian categories, thus providing an abstract framework for retrodiction beyond classical and quantum probability theory. 
More precisely, we extend Bayesian inference, and more generally Jeffrey's probability kinematics, to arbitrary semicartesian categories. 


\end{abstract}

\vspace{-7mm}
\tableofcontents

\section{Introduction and main results}

Categorical probability theory aims to isolate categorical~\cite{Ma98} axioms that can be used to reason about classical information flow, statistics, probability theory, entropy, etc.\ through a set of logical, synthetic, and oftentimes operational, axioms~\cite{La62,Gi82,Ga96,Fo12,BaFr14,CuSt14,DSDG18,ChJa18,Fr20}. 
Similarly, process theories~\cite{AbCo04,Chiribella14distinguishability,SSC21,DDR21}, which are modeled by monoidal categories~\cite{Ya24,JSV96,Se10,Ma98}, are used as general models for physical theories of information flow that are classical, quantum, and beyond these two regimes in order to better understand possible generalizations of information theory itself. 
Special examples include general probabilistic theories (GPTs)~\cite{Hardy01,CDP10} (see the introduction in Ref.~\cite{Chiribella14distinguishability} for a list of relevant references). 
These perspectives are somewhat analogous to, yet quite distinct from, the development of non-commutative probability through algebraic reformulations of concepts from probability theory, statistics, and information theory~\cite{Mo54,Um54}. 
The similarities between the two perspectives of non-commutative probability and categorical probability have recently been discussed in Ref.~\cite{PaBayes,GPRR21}. 

Many classical definitions and foundational results from probability and statistics have been extended to and derived from such categorical axioms, including Shannon entropy~\cite{BFL}, relative entropy (Kullback--Liebler divergence)~\cite{BaFr14,GaPa18}, conditional entropy~\cite{FuPa21}, mutual information~\cite{Fu23}, statistical divergences~\cite{Perrone23}, conditioning~\cite{ChJa18,Fr20}, independence~\cite{ChJa18,Fr20}, disintegrations~\cite{PaRu19,PaBayes}, Bayesian networks and causal inference~\cite{JKZ19,JaZa20,FrKl23}, de~Finetti's theorem~\cite{FGP21,MoPe22,JaSt20}, theorems on sufficient statistics~\cite{Fr20}, the Kolmogorov and Hewitt--Savage zero-one laws~\cite{FrRi20}, statistical comparison~\cite{FGPR23}, the $d$-separation criterion~\cite{FrKl23}, the Ergodic Decomposition Theorem~\cite{MoPe23}, and more. This plethora of results illustrates that such classical results and their proofs can occasionally be obtained by synthetic means rather than choosing specific analytical models for probability (see the introduction of Ref.~\cite{MoPe23} for an illuminating discussion about synthetic versus analytic). Categorical perspectives on probability theory and information theory have also been helpful in computer science for finding more adequate theories of randomness in function spaces~\cite{HKSY17,SSSW21}. In addition, some of the aforementioned results have been extended to the quantum setting as well~\cite{PaEntropy,PaRelEnt,PaBayes,PaQPL21,StSu23,FrLo23}. These illustrate some of the benefits of a categorical approach, namely the discovery of new results, methods, and techniques in classical or quantum probability and/or information theory. 

However, more work needs to be done in order to further develop the theory. In particular, the notion of Bayesian inversion as introduced in the categorical probability literature, specifically Ref.~\cite{ChJa18} (see also Refs.~\cite{Fo12,CuSt14,Fr20,DDGK16,CDDG17,DSDG18}), does not extend in a fully robust manner to quantum theories, specifically quantum channels, which model stochastic evolution in quantum systems~\cite{Kr83,NiCh11}. Namely, Refs.~\cite{PaRuBayes,GPRR21,PaRu19} show that additional symmetry constraints are required for a state-preserving quantum channel to have a Bayesian inverse with respect to the categorical definition from Ref.~\cite{ChJa18}. These symmetry constraints are related to the notion of time-covariance of quantum channels~\cite{Ma12}, or equivalently invariance under the modular Hamiltonian flow associated with the initial and final states~\cite{AcCe82}, which is a generalization of GNS symmetry coming from the GNS inner product detailed balance condition~\cite{FaUm07} (which is distinct from the KMS and Bogoliubov inner products~\cite{Ha17,HHW67,PeTo93,CaMa19}). Despite these negative results, many justifications for extensions of Bayes' rule and Bayesian inversion, sometimes called \emph{time-reversal symmetry} or \emph{retrodiction} due to their usage in making inference about the past, have been put forward in the context of \emph{recovery maps}~\cite{BaKn02,Cr08,Le06,Le07,Ty10,LeSp13,LiWi18,SuToHa16,CHPSSW19,JRSWW16,JRSWW18,FSB20,SuSc22,BuSc21,AwBuSc21,Ts22}. Many, though not all of these, are based off of the Petz recovery map or its variants~\cite{Pe84,Wilde15}. When the aforementioned symmetry condition holds, many of these proposals coincide~\cite{GPRR21,PaBu22}. Most of these proposed versions of Bayes' rule have now been shown to be special cases of a fully quantum Bayes' rule~\cite{FuPa22a} that depends on a quantum extension of joint probabilities, and which are called \emph{states over time}~\cite{HHPBS17,FuPa22} (see also Refs.~\cite{Ts22,Ts23,PFBC23,LiNg23} for recent developments and alternative sources).

However, in Ref.~\cite{PaBu22}, it was shown that only \emph{one} of the recovery maps for quantum channels from the previous paragraph satisfies dagger monoidal functoriality~\cite{Ka18} properties of classical Bayesian inversion. These functoriality properties were proved in Ref.~\cite{DSDG18} in the context of classical probability and statistics and a subset of these properties were proved independently in Refs.~\cite{LiWi18,Wilde15} in the context of quantum information theory for their desirable properties for approximate quantum error-correction. 
These functoriality properties have also been shown to be important for understanding multipartite entanglement in conformal field theories involving multiple regions~\cite{VWZ23}. 
Later, Ref.~\cite{Fr20} proved a rather general theorem that says any (classical) Markov category with conditionals admits a dagger on its associated category of states and processes, which specializes to the case of Bayesian inversion in the context of classical probability theory (see also Ref.~\cite{EnPe23}) and quantum systems with symmetry (which includes classical systems)~\cite{PaBayes}, but not to all (general) quantum systems~\cite{PaRuBayes,GPRR21,PaBayes,PaRu19,PaBu22}. 

The purpose of the present paper is to form a rigorous categorical, i.e., process-theoretic, foundation for retrodiction using semicartesian categories as models for information flow~\cite{FGGPS22} and dagger categories~\cite{Ka18} as models for the reversal of such information flow in a way that specializes to the known cases of classical Bayesian inversion and quantum retrodiction through the Petz recovery map~\cite{Pe84,PaBu22}. Our definition can be viewed as a way to reverse information flow in such models, though it is more accurate to view this reversal as retrodiction or backwards belief propagation rather than physical reversals of evolution~\cite{Wat55,LeSp13,PaBayes,BuSc21,AwBuSc21,PaBu22,Cr08}. 

The paper is organized as follows. Section~\ref{sec:classicalBayes} goes over the prototypical example of a semicartesian category, namely stochastic maps, i.e., conditional probabilities, between finite sets. It also reviews Bayesian inversion and describes the main categorical properties that Bayesian inversion satisfies. Section~\ref{sec:semicartesiancats} reviews semicartesian categories, dilations, dilational equality, and faithful states, and it proves several new results concerning the semicartesian category of finite-dimensional $C^*$-algebras and completely positive trace-preserving maps, i.e., quantum channels (the usage of $C^*$-algebras, as opposed to just Hilbert spaces, is to allow for hybrid classical/quantum systems and channels, which includes preparations, measurements, and instruments).   Section~\ref{sec:retrodiction} defines retrodiction and probability kinematics abstractly in the context of semicartesian categories. Whereas Ref.~\cite{PaBu22} proved the existence of a retrodiction functor on a category of \emph{faithful} states and state-preserving quantum channels, Section~\ref{sec:retrodiction} proves the existence of a retrodiction functor on the category of \emph{all} states and dilationally equal equivalence classes of state-preserving quantum channels.

\section{Bayesian retrodiction in classical probability}
\label{sec:classicalBayes}

We begin our study of Bayesian retrodiction by first recalling a well-established case in the context of finite probability theory and statistics~\cite{Ba1763}. However, our formulation of Bayesian retrodiction and Jeffrey's probability kinematics~\cite{Je90,Ja19,Ja20}, which describes how probabilities are updated based on probabilistic evidence, are emphasized through their categorical and compositional properties~\cite{PaBu22}.

\bd
Given two finite sets $X$ and $Y$, a \define{stochastic map} from $X$ to $Y$, written as $X\xrightarrow{f}Y$ or $f:X\to Y$, is an assignment sending each $x\in X$ to a probability distribution $f_{x}$ on $Y$ whose value on $y\in Y$ is denoted by $f_{yx}$. Thus, $f_{yx}\ge0$ for all $x\in X$ and $y\in Y$ and $\sum_{y\in Y}f_{yx}=1$ for all $x\in X$. Denoting a single element set by $\{\bullet\}$, a stochastic map $\{\bullet\}\xrightarrow{p}X$ assigns to the only element $\bullet$ a probability distribution on $X$, i.e., $p$ is just a single probability distribution on $X$. As such, its value on $x$ is denoted by $p_{x}$ (as opposed to $p_{x\bullet}$).
\ed

\br
If $X\xrightarrow{f}Y$ is a \emph{function} between finite sets, it defines a stochastic map $X\xrightarrow{\delta(g)}Y$ upon setting $\delta(g)_{yx}:=\delta_{yg(x)}$. By a slight abuse of notation, $\delta(g)$ will often be denoted by $g$ so that functions will be viewed as stochastic maps upon this identification. 
\er

\bd
Given finite sets $X,Y,$ and $Z$, and a pair of successive stochastic maps $X\xrightarrow{f}Y\xrightarrow{g}Z$, their \define{composite} is the stochastic map $g\circ f$ given by matrix multiplication, i.e., 
\[
(g\circ f)_{zx}:=\sum_{y\in Y}g_{zy}f_{yx},
\]
which is called the \emph{Chapman--Kolmogorov equation}. In particular, if $\{\bullet\}\xrightarrow{p}X$ is a probability distribution on $X$, then $q:=f\circ p$ is the \define{pushforward} probability distribution on $Y$, defined by the formula $q_{y}=\sum_{x\in X}f_{yx}p_{x}$ for all $y\in Y$. 
\ed

\bd
Let $\FinStoch$ be the symmetric monoidal category~\cite{Ya24,Se10,Ma98} whose objects are finite sets, whose morphisms are stochastic maps, and where the composition is defined by the Chapman--Kolmogorov equation~\cite{BaFr14}. The monoidal structure on objects is given by $X\otimes X':=X\times X'$, the Cartesian product. Meanwhile, the monoidal product on morphisms $X\xrightarrow{f}Y$ and $X'\xrightarrow{f'}Y'$ is given by the stochastic map $X\times X'\xrightarrow{f\otimes f'}Y\times Y'$ whose components are
\[
(f\otimes f')_{(y,y')(x,x')}:=f_{yx}f_{y'x'}
\]
for all $x\in X, x'\in X', y\in Y,$ and $y'\in Y'$. The braiding isomorphism $X\times Y\xrightarrow{\gamma}Y\times X$, also called the \define{swap} isomorphism%
\footnote{This is not to be confused with the ${\tt SWAP}$ gate of quantum circuits.}
 is the standard one corresponding to the function sending $(x,y)$ to $(y,x)$. The set $\{\bullet\}$ is a monoidal unit. 
\ed

\bd
Let $\stFS$ be the symmetric monoidal category whose objects are pairs $(X,p)$ with $X$ a finite set and $p$ a faithful (nowhere vanishing) probability distribution on $X$, i.e., $p_{x}>0$ for all $x\in X$. A morphism $f$ from $(X,p)$ to $(Y,q)$ in $\stFS$ is a stochastic map%
\footnote{In Refs.~\cite{BaFr14,PaRelEnt}, a similar, but different, category $\FinProb$ is defined. Namely, the objects of $\FinProb$ are still pairs $(X,p)$, but the probabilities $p$ are arbitrary, i.e., $p_{x}\ge0$ for all $x\in X$. Second, the morphisms  $(X,p)\xrightarrow{f}(Y,q)$ of $\FinProb$ are required to be \define{deterministic} in the sense that $f_{yx}\in\{0,1\}$ for all $x\in X$ and $y\in Y$.}
 $X\xrightarrow{f}Y$ such that $f\circ p=q$, sometimes called a \define{state-preserving morphism} to emphasize this condition. 
The composition in $\stFS$ is given by the composition of the underlying stochastic maps.
The tensor product $(X,p)\otimes (X',p')$ of objects is given by the cartesian product for the underlying sets and the product of probability distributions $(X,p)\otimes (X',p'):=(X\times X',p\otimes p')$, where $(p\otimes p')_{(x,x')}:=p_{x}p'_{x'}$ for all $x\in X$ and $x'\in X'$. 
Meanwhile, the tensor product of $(X,p)\xrightarrow{f}(Y,q)$ and $(X',p')\xrightarrow{f'}(Y',q')$ is given by the morphism $(X\times X',p\otimes p')\xrightarrow{f\otimes f'}(Y\times Y',q\otimes q')$ (it follows from this definition that the \define{interchange law} $(f\otimes f')\circ(p\otimes p')=q\otimes q'$ holds). 
\ed

\br
Note that the definition of $\stFS$ is not exactly the coslice/under category~\cite{Ma98} $\{\bullet\}\downarrow\FinStoch$ because the probability distributions $p$ from the objects $(X,p)$ are required to be nowhere vanishing. As such, $\stFS$ is a \emph{subcategory} of this coslice category. One can therefore visualize $(X,p)$ and $(X,p)\xrightarrow{f}(Y,q)$ as
\[
\xy0;/r.25pc/:
(0,6)*+{\{\bullet\}}="1";
(0,-6)*+{X}="X";
{\ar"1";"X"_{p}};
\endxy
\qquad\text{ and }\qquad
\xy0;/r.25pc/:
(0,6)*+{\{\bullet\}}="1";
(-6,-6)*+{X}="X";
(6,-6)*+{Y}="Y";
{\ar"1";"X"_{p}};
{\ar"1";"Y"^{q}};
{\ar"X";"Y"_{f}};
\endxy
\;\;,
\]
respectively, where the diagram on the right commutes in the sense that the stochastic map $f\circ p$ equals the stochastic map $q$. 
\er

\bd
Given a morphism $(X,p)\xrightarrow{f}(Y,q)$ in $\stFS$, the probability distribution $p$ is called the \define{prior}, the stochastic map $f$ is called the \define{likelihood}, and $q$ is called the \define{prediction} or \define{marginal likelihood}. Given such a morphism $f$, its \define{Bayesian inverse} is the morphism $(X,p)\xleftarrow{\ov f}(Y,q)$ in $\stFS$ whose components are given by the formula 
\be
\label{eq:Bayesclassical}
\ov f_{xy}:=\frac{f_{yx}p_{x}}{q_{y}}.
\ee
\ed

\br
One of the reasons $\ov f$ is called a Bayesian \emph{inverse} is due to the fact that if $(X,p)\xrightarrow{f}(Y,q)$ is an invertible morphism in $\stFS$, then $\ov f=f^{-1}$, where $(X,p)\xleftarrow{f^{-1}}(Y,q)$ is the inverse of $f$. 
This shows that the \emph{procedure} of constructing a Bayesian inverse $f\mapsto\ov f$ is a generalization of the \emph{procedure} of constructing the inverse of an invertible morphism $f\mapsto f^{-1}$. In fact, the common properties of these procedures are used to \emph{define} retrodiction later in this paper. Importantly, although not every morphism in $\stFS$ is invertible, every morphism is \emph{Bayesian invertible}. 
\er

\begin{term}
Given a morphism $(X,p)\xrightarrow{f}(Y,q)$ in $\stFS$, if $e$ is another probability distribution on $Y$, called \define{soft evidence} (and interpreted as an observation), then 
\[
(\ov f\circ e)_{x}=\sum_{y\in Y}\frac{f_{yx}p_{x}e_{y}}{q_{y}}
\]
is the \define{updated} probability distribution based on the evidence $e$. This procedure of pushing evidence $e$ backwards along $\ov f$ is called \define{Jeffrey's probability kinematics} and generalizes the standard Bayes' rule, the latter of which deals with the special case of \define{hard evidence} given by a \define{Dirac delta} distribution $e=\delta_{y}$ for some $y\in Y$ and which coincides with~\eqref{eq:Bayesclassical}~\cite{Je90,Ja19,Ja20}. One can visualize Jeffrey's probability kinematics schematically via 
\[
\xy0;/r.25pc/:
(0,6)*+{\{\bullet\}}="1";
(-6,-6)*+{X}="X";
(6,-6)*+{Y}="Y";
{\ar"1";"X"_{p}};
{\ar"1";"Y"^{q}};
{\ar"X";"Y"_{f}};
\endxy
\quad\text{ and }\quad
\xy0;/r.25pc/:
(0,6)*+{\{\bullet\}}="1";
(6,-6)*+{Y}="Y";
{\ar"1";"Y"^{e}};
\endxy
\quad\text{ give }\quad
\xy0;/r.25pc/:
(0,6)*+{\{\bullet\}}="1";
(-6,-6)*+{X}="X";
(6,-6)*+{Y}="Y";
{\ar"1";"X"_{\ov f\circ e}};
{\ar"1";"Y"^{e}};
{\ar"Y";"X"^{\ov f}};
\endxy
\]
\end{term}

\bt
\label{thm:BisRetro}
The assignment 
$\stFS\xrightarrow{\retro^{\Bayes}}\stFS^{\op}$ sending a morphism to its Bayesian inverse 
satisfies the following properties: 
\begin{enumerate}
\item
\label{item:recovery}
(recovery property) $\retro^{\Bayes}$ acts as the identity on objects (and thus sends a morphism $(X,p)\xrightarrow{f}(Y,q)$ to another morphism $(X,p)\xleftarrow{\retro^{\Bayes}(f)}(Y,q)$). 
\item
(identity-preservation) $\retro^{\Bayes}$ fixes each identity morphism, i.e., $\retro^{\Bayes}(\id_{(X,p)})=\id_{(X,p)}$ for all objects $(X,p)$. 
\item
(compositionality) For every pair $(X,p)\xrightarrow{f}(Y,q)\xrightarrow{g}(Z,r)$ of composable morphisms, $\retro^{\Bayes}(g\circ f)=\retro^{\Bayes}(f)\circ\retro^{\Bayes}(g)$. 
\item
(tensoriality) $\retro^{\Bayes}(f\otimes f')=\retro^{\Bayes}(f)\otimes\retro^{\Bayes}(f')$ for every pair of morphisms $(X,p)\xrightarrow{f}(Y,q)$ and $(X',p')\xrightarrow{f'}(Y',q')$. 
\item
(extending inversion) $\retro^{\Bayes}$ is \define{inverting}, i.e., $\retro^{\Bayes}(f)=f^{-1}$ whenever $f$ is an isomorphism%
\footnote{Identify-preservation is a special case of this axiom.}.
\item
\label{item:involutivity}
(involutivity) $\retro^{\Bayes}$ is \define{involutive}, i.e., $\retro^{\Bayes}\circ\retro^{\Bayes}=\id_{\stFS}$.
\end{enumerate}
In other words, $\stFS$ equipped with $\retro^{\Bayes}$ is a symmetric monoidal inverting dagger category~\cite{Ka18,PaBu22}.  
\et

According to Ref.~\cite{HeKo22}, Ref.~\cite{Ma50} first introduced dagger categories under a different name and Ref.~\cite{Se07} then introduced the currently used terminology. Ref.~\cite{PaBu22} emphasized the importance of including the inverting condition. Ref.~\cite{Ka18} is a thorough monograph on dagger categories. 

\br
Identity-preservation and compositionality say that $\retro^{\Bayes}$ is a functor (the ${}^{\op}$ says it's technically a contravariant functor~\cite{Ma98}). The recovery property additionally says that $\retro^{\Bayes}$ is an identity-on-objects functor. Tensoriality and functoriality say that $\retro^{\Bayes}$ is a monoidal functor. Including the recovery and involutivity properties additionally say that $\retro^{\Bayes}$ is a monoidal dagger (sometimes called a monoidal dagger functor). The inverting condition is one of the key axioms indicating how $\retro^{\Bayes}$ provides an extension of inversion to all morphisms. This is stronger than the category $\stFS$ being unitary, in the terminology of Ref.~\cite{Vi11}, which is a condition that says if two objects $(X,p)$ and $(Y,q)$ are isomorphic, then there \emph{exists} a morphism $(X,p)\xrightarrow{f}(Y,q)$ such that $\retro^{\Bayes}(f)=f^{-1}$. The inverting condition says that \emph{every} isomorphism is unitary in the sense of Ref.~\cite{Vi11}. 
\er

\bd
An assignment $\stFS\xrightarrow{\retro}\stFS^{\op}$ satisfying properties~\ref{item:recovery}--\ref{item:involutivity} in Theorem~\ref{thm:BisRetro} is called a \define{retrodiction functor} on $\stFS$. 
\ed

Versions of Theorem~\ref{thm:BisRetro} have appeared in many forms throughout the literature such as~\cite[Remark~13.10]{Fr20} and~\cite[Theorem 2.10]{DSDG18} (see Example~\ref{ex:FinStochretro} for more details), though our present formulation is based on Ref.~\cite{PaBu22}, emphasizing its role as an extension of time-reversal symmetry.  
A special case of the main open question asked in Ref.~\cite{PaBu22} is the following.

\bq
\label{qu:Bunique}
Given a retrodiction functor $\stFS\xrightarrow{\retro}\stFS^{\op}$, does $\retro=\retro^{\Bayes}$?
\eq

Theorem~\ref{thm:BisRetro} tells us that Bayesian inversion gives \emph{one example} of a retrodiction functor $\stFS\to\stFS^{\op}$. Meanwhile, Question~\ref{qu:Bunique} is asking if Bayesian inversion is the \emph{only} retrodiction functor $\stFS\to\stFS^{\op}$. Thus, the question of characterizing classical Bayesian inference in terms of logical axioms has now been formulated as a precise mathematical question. Note that it is different from earlier characterizations of Bayesian inference~\cite{Cs91}, and a new (unfinished) attempt in terms of detailed balance and a certain optimality condition involving matrix determinants~\cite{SuSc22}. In particular, Question~\ref{qu:Bunique} is an attempt at a purely categorical characterization of Bayesian inference without the need for a specific model and without optimizing a specific quantity (whose choice is somewhat arbitrary). It therefore completely does away with probabilities in its description and is thus a purely process-theoretic characterization in the spirit of Refs.~\cite{Chiribella14distinguishability,Fr20} going beyond general probabilistic theories (GPTs).

\section{Semicartesian categories and state-preserving channels}
\label{sec:semicartesiancats}

In what follows, we will model physical theories involving states and evolution by semicartesian categories, primarily referring to Refs.~\cite{JSV96,Se10,FGGPS22,Ho21,Ya24} as guides. Briefly, semicartesian categories are symmetric monoidal categories where the monoidal unit is terminal (this is sometimes called a \emph{causal} process theory in the literature~\cite{CDP10,SSC21}). They can be used to model both classical and quantum systems, including stochastic aspects, such as density matrices, preparations, measurements, quantum channels, and instruments~\cite{SSGC22,FuPa22a}. 

\bd
A \define{semicartesian category} is a symmetric monoidal category $\mathbf{C}$ whose monoidal unit $\Alg{I}$ is terminal. The morphism from $\Alg{A}$ to $\Alg{I}$ is denoted by $!_{\Alg{A}}$ and is called the \define{grounding} morphism. A \define{state} in a semicartesian category is a morphism whose domain is $\Alg{I}$. 
Given a finite set $\mathcal{I}$, the map $\bigotimes_{i\in\mathcal{I}}\Alg{A}_{i}$ to $\bigotimes_{j\in\mathcal{J}}\Alg{A}_{j}$ with $\mathcal{J}\subseteq\mathcal{I}$ given by grounding all $\Alg{A}_{k}$ with $k\notin\mathcal{K}:=\mathcal{I}\setminus\mathcal{J}$ is denoted by $!_{\bigotimes_{k\in\mathcal{K}}\Alg{A}_{k}}$ and is called the \define{partial grounding}. 
\ed

For example, if $\mathcal{I}=\{1,2,3,4,5\}$ and $\mathcal{J}=\{1,3,4\}$, then $\mathcal{K}=\{2,5\}$ and 
\[
!_{\Alg{A}_{2}\otimes\Alg{A}_{5}}=\id_{\Alg{A}_{1}}\otimes!_{\Alg{A}_{2}}\otimes\id_{\Alg{A}_{3}}\otimes\id_{\Alg{A}_{4}}\otimes!_{\Alg{A}_{5}}.
\]
String-diagram visualizations of some morphisms in a semicartesian category are shown in Table~\ref{table:stringdiagrams}. The reader is referred to Refs.~\cite{JSV96,Se10,ChJa18,FGGPS22,SSGC22} for additional details on string diagrams as used in this paper. Nevertheless, the reader need not be familiar with them to understand the main claims made in this paper. 

\begin{table}
\centering
\begin{tabular}{ccccccc}
$\vcenter{\hbox{%
\begin{tikzpicture}[font=\small]
\node[arrow box] (c) at (0,0) {$\map{E}$};
\draw (c) to (0,0.8);
\draw (c) to (0,-0.8);
\node at (-0.2,-0.70) {\scriptsize $\Alg{A}$};
\node at (-0.2,0.70) {\scriptsize $\Alg{B}$};
\end{tikzpicture}}}$
&
$\vcenter{\hbox{%
\begin{tikzpicture}[font=\small]
\draw (0,-0.7) to (0,0.7);
\node at (-0.2,0) {\scriptsize $\Alg{A}$};
\end{tikzpicture}}}$
&
$\vcenter{\hbox{%
\begin{tikzpicture}[font=\small]
\node[discarder] (c) at (0,0) {};
\draw (c) to (0,-0.8);
\node at (-0.2,-0.70) {\scriptsize $\Alg{A}$};
\end{tikzpicture}}}$
&
$\vcenter{\hbox{%
\begin{tikzpicture}[font=\small]
\node[state] (omega) at (0,-0.5) {$\alpha$};
\coordinate (q) at (0,0.2);
\draw (omega) to (q);
\node at (-0.2,0.10) {\scriptsize $\Alg{A}$};
\end{tikzpicture}}}$
&
$\vcenter{\hbox{%
\begin{tikzpicture}[font=\small]
\node[arrow box] (c) at (0,0) {$\map{E}$};
\draw (c) to (0,0.8);
\draw (c) to (0,-0.8);
\node at (-0.2,-0.70) {\scriptsize $\Alg{A}$};
\node at (-0.2,0.70) {\scriptsize $\Alg{B}$};
\node[arrow box] (c2) at (0.8,0) {$\map{E}'$};
\draw (c2) to (0.8,0.8);
\draw (c2) to (0.8,-0.8);
\node at (0.6,-0.70) {\scriptsize $\Alg{A}'$};
\node at (0.6,0.70) {\scriptsize $\Alg{B}'$};
\end{tikzpicture}}}$
&
$\vcenter{\hbox{%
\begin{tikzpicture}[font=\small]
\coordinate (X) at (-0.45,0) {};
\coordinate (Y) at (0.45,0) {};
\coordinate (X2) at (0.45,1.2) {};
\coordinate (Y2) at (-0.45,1.2) {};
\draw (X) to[out=90,in=-90] (X2);
\draw (Y) to[out=90,in=-90] (Y2);
\path[scriptstyle]
node at (-0.65,0.15) {$\Alg{A}$}
node at (0.65,0.15) {$\Alg{B}$}
node at (-0.65,1.05) {$\Alg{B}$}
node at (0.65,1.05) {$\Alg{A}$};
\end{tikzpicture}}}$
&
$\vcenter{\hbox{%
\begin{tikzpicture}[font=\small]
\draw (0,-0.7) to (0,0.7);
\node at (-0.2,-0.5) {\scriptsize $\Alg{A}_{1}$};
\node[discarder] (2) at (0.55,0.1) {};
\draw (0.55,-0.7) to (2);
\node at (0.35,-0.5) {\scriptsize $\Alg{A}_{2}$};
\draw (1.1,-0.7) to (1.1,0.7);
\node at (0.9,-0.5) {\scriptsize $\Alg{A}_{3}$};
\draw (1.65,-0.7) to (1.65,0.7);
\node at (1.45,-0.5) {\scriptsize $\Alg{A}_{4}$};
\node[discarder] (5) at (2.2,0.1) {};
\draw (2.2,-0.7) to (5);
\node at (2.0,-0.5) {\scriptsize $\Alg{A}_{5}$};
\end{tikzpicture}}}$
\\
\begin{tabular}{c}arbitrary\\morphism\end{tabular}&identity&grounding&state&\begin{tabular}{c}tensor\\product\end{tabular}&\begin{tabular}{c}swapping/\\braiding\end{tabular}&\begin{tabular}{c}partial\\grounding\end{tabular}
\\
$\Alg{A}\xrightarrow{\map{E}}\Alg{B}$
&
$\Alg{A}\xrightarrow{\id_{\Alg{A}}}\Alg{A}$
&
$\Alg{A}\xrightarrow{!_{\Alg{A}}}I$
&
$I\xrightarrow{\alpha}\Alg{A}$
&
$\map{E}\otimes\map{E}'$
&
$\Alg{A}\otimes\Alg{B}\xrightarrow{\gamma}\Alg{B}\otimes\Alg{A}$
&
$!_{\Alg{A}_{2}\otimes\Alg{A}_{5}}$
\end{tabular}
\caption{Several examples of morphisms in a semicartesian category along with string-diagrammatic representations. The partial grounding is a morphism of the form $\Alg{A}_{1}\otimes\Alg{A}_{2}\otimes\Alg{A}_{3}\otimes\Alg{A}_{4}\otimes\Alg{A}_{5}\xrightarrow{!_{\Alg{A}_{2}\otimes\Alg{A}_{5}}}\Alg{A}_{1}\otimes\Alg{A}_{3}\otimes\Alg{A}_{4}$}
\label{table:stringdiagrams}
\end{table}

\bx
$\FinStoch$ is a semicartesian category (see Ref.~\cite{FGGPS22} for details).  
\ex

\bx
Let $\NCFinStoch$ be the category whose objects are finite-dimensional unital $C^*$-algebras~\cite{Fa01} (as all $C^*$-algebras here will be unital, we henceforth drop the word unital) and whose morphisms are completely positive trace-preserving (CPTP) maps~\cite{Ch75,Kr83}. Composition is defined by function composition and the tensor product is the usual one for algebras (by the finite-dimensionality assumption, the norm on the tensor product is canonically defined~\cite{Fa01}). 
With these specifications, $\NCFinStoch$ is a semicartesian category. 
The grounding morphism is the trace, while the partial grounding is the partial trace. 
As such, this grounding is denoted by $\Tr$ instead of $!$. 
\ex

\br
Note that there is a fully faithful functor $\FinStoch\hookrightarrow\NCFinStoch$ given by sending $X$ to $\C^{X}$ and sending $X\xrightarrow{f}Y$ to the CPTP map $\C^{X}\xrightarrow{\map{F}}\C^{Y}$ uniquely determined by sending $\delta_{x}\in\C^{X}$, the Kronecker delta function at $x\in X$, the value of which is $1$ at $x$ and $0$ otherwise, to the function $\map{F}(\delta_{x})\in\C^{Y}$ on $Y$, the value of which is given by $\map{F}(\delta_{x})(y):=f_{yx}$ for $y\in Y$~\cite{Pa17}.  
Note that our definition of $\NCFinStoch$ differs from the definition of the semicartesian category $\mathbf{QIT}$ of Ref.~\cite{Ho21}, the latter of which uses finite-dimensional Hilbert spaces as objects and quantum channels (CPTP maps) between the $C^*$-algebras of bounded operators on those Hilbert spaces as the morphisms. This causes some issues in terms of embedding classical systems into quantum systems (cf.\ \cite[Remark~1.1.12]{Ho21}) because the embedding is not a functor. By working with \emph{algebras} instead of \emph{Hilbert spaces}, we are working in a category of \emph{hybrid} classical/quantum systems, and we therefore naturally obtain a fully faithful functor from classical information theory into quantum information theory that gives an equivalence on the subcategory of $\NCFinStoch$ consisting of \emph{commutative} $C^*$-algebras. This equivalence is sometimes called \emph{probabilistic Gelfand duality} or the \emph{stochastic commutative Gelfand--Naimark theorem}~\cite{Pa17,FuJa15}. Thus, classical systems embed functorially into classical/quantum systems provided that one uses the category $\NCFinStoch$, which is what we do in this paper.
\er

We will soon define categories of states and state-preserving processes. But before we do that, we must say what it means for a state to be faithful. This can be done using the theory of dilations by an appropriate notion of almost surely (a.s.) equivalence~\cite{CDP10,Chiribella14dilation,SSC21,Ho21,FGGPS22}. 

\bd
Let $\mathbf{C}$ be a semicartesian category and let $\Alg{A}\xrightarrow{\map{E}}\Alg{B}$ be a morphism in $\mathbf{C}$. A \define{dilation} of $\map{E}$ is a morphism of the form $\Alg{A}\xrightarrow{\map{D}}\Alg{B}\otimes\Alg{E}$ 
satisfying 
$!_{\Alg{E}}\circ\map{D}=\map{E}$, i.e., 
\[
\vcenter{\hbox{
\begin{tikzpicture}[font=\small]
\node[arrow box] (p) at (0,0) {\;\;$\map{D}$\;\;};
\coordinate (X) at (-0.25,0.8);
\coordinate (Y) at (0.25,0.8);
\draw (0,-0.8) to (p);
\draw (p.north)++(-0.25,0) to (X);
\draw (p.north)++(0.25,0) to (Y);
\node at (-0.2,-0.70) {\scriptsize $\Alg{A}$};
\node at (-0.4,0.70) {\scriptsize $\Alg{B}$};
\node at (0.4,0.70) {\scriptsize $\Alg{E}$};
\end{tikzpicture}}}
\quad\text{ such that }\quad
\vcenter{\hbox{
\begin{tikzpicture}[font=\small]
\node[arrow box] (p) at (0,0) {\;\;$\map{D}$\;\;};
\coordinate (X) at (-0.25,0.8);
\node[discarder] (Y) at (0.25,0.5) {};
\draw (0,-0.8) to (p);
\draw (p.north)++(-0.25,0) to (X);
\draw (p.north)++(0.25,0) to (Y);
\node at (-0.2,-0.70) {\scriptsize $\Alg{A}$};
\node at (-0.4,0.70) {\scriptsize $\Alg{B}$};
\end{tikzpicture}}}
\;\;=\;\;
\vcenter{\hbox{%
\begin{tikzpicture}[font=\small]
\node[arrow box] (c) at (0,0) {$\map{E}$};
\draw (c) to (0,0.8);
\draw (c) to (0,-0.8);
\node at (-0.2,-0.70) {\scriptsize $\Alg{A}$};
\node at (-0.2,0.70) {\scriptsize $\Alg{B}$};
\end{tikzpicture}}}
\;\;.
\]
The object $\Alg{E}$ will be called the \define{environment} associated with the dilation $\map{D}$. 
\ed

\br
The terminology of dilation here follows the one of Refs.~\cite{SSC21,FGGPS22}, while Ref.~\cite{Ho21} calls this a \emph{one-sided dilation}. Since one-sided dilations are the only ones we will work with here, we prefer to use the more concise terminology. 
\er

Before giving examples of dilations in $\FinStoch$ and $\NCFinStoch$, we define pure states. Purity of processes, and more specifically states, can also be defined process-theoretically~\cite{Chiribella14distinguishability,SSC21}. In this paper, we will only use the definition of a pure state, rather than for arbitrary processes. 

\bd
A state $I\xrightarrow{\alpha}\Alg{A}$ in a semicartesian category $\mathbf{C}$ is \define{pure} iff for any dilation $I\xrightarrow{\pi}\Alg{A}\otimes\Alg{E}$ of $\alpha$, there exists a state $I\xrightarrow{\epsilon}\Alg{E}$ such that $\pi=\alpha\otimes\epsilon$, i.e., 
\[
\vcenter{\hbox{%
\begin{tikzpicture}[font=\small]
\node[state] (omega) at (0,0) {\;$\pi$\;};
\coordinate (X) at (-0.25,0.95) {};
\node[discarder] (Y) at (0.25,0.65) {};
\draw (omega) ++(-0.25, 0) to (X);
\draw (omega) ++(0.25, 0) to (Y);
\path[scriptstyle]
node at (-0.45,0.75) {$\Alg{A}$}
node at (0.45,0.35) {$\Alg{E}$};
\end{tikzpicture}}}
\;\;=\;\;
\vcenter{\hbox{%
\begin{tikzpicture}[font=\small]
\node[state] (omega) at (0,-0.8) {$\alpha$};
\coordinate (q) at (0,-0.1);
\draw (omega) to (q);
\path[scriptstyle]
node at (-0.20,-0.25) {$\Alg{A}$};
\end{tikzpicture}}}
\qquad\implies\qquad
\vcenter{\hbox{%
\begin{tikzpicture}[font=\small]
\node[state] (omega) at (0,0) {\;$\pi$\;};
\coordinate (X) at (-0.25,0.95) {};
\coordinate (Y) at (0.25,0.95) {};
\draw (omega) ++(-0.25, 0) to (X);
\draw (omega) ++(0.25, 0) to (Y);
\path[scriptstyle]
node at (-0.45,0.75) {$\Alg{A}$}
node at (0.45,0.75) {$\Alg{E}$};
\end{tikzpicture}}}
\;\;=\;\;
\vcenter{\hbox{%
\begin{tikzpicture}[font=\small]
\node[state] (omega) at (0,-0.8) {$\alpha$};
\coordinate (q) at (0,-0.1);
\node[state] (xi) at (1.0,-0.8) {$\epsilon$};
\coordinate (p) at (1.0,-0.1);
\draw (omega) to (q);
\draw (xi) to (p);
\path[scriptstyle]
node at (-0.20,-0.25) {$\Alg{A}$}
node at (0.80,-0.25) {$\Alg{E}$};
\end{tikzpicture}}}
\]
for some state $I\xrightarrow{\epsilon}\Alg{E}$. A \define{purification} of $\alpha$ is a dilation $I\xrightarrow{\pi}\Alg{A}\otimes\Alg{E}$ of $\alpha$ with $\pi$ a pure state. 
\ed

\bx
A pure state in $\FinStoch$ is a Dirac delta probability measure, while a pure state in $\NCFinStoch$ is a pure state in the sense of non-commutative probability, namely an extreme point of the convex space of all states~\cite{BeZy06,Chiribella14distinguishability,SSC21}.
\ex

\bx
\label{ex:dilationsCAlg}
Let $\Alg{A}$ be an object of $\NCFinStoch$, i.e., a finite-dimensional $C^*$-algebra, and let $\alpha\in\Alg{A}$ be an arbitrary state, specified as the image of $1\in\C$ under a CPTP map $\C\xrightarrow{\alpha}\Alg{A}$, which is somewhat abusively also be denoted by $\alpha$. We will look at four examples of dilations. 
\begin{enumerate}[(a)]
\item
\label{item:purification}
Say $\Alg{A}=\matr_{m}$, the $C^*$-algebra of $m\times m$ matrices. Write $\alpha=\sum_{i=1}^{m}p_{i}|i\>\<i|$ in terms of an orthonormal basis of eigenvectors of $\alpha$ so that the $\{p_{i}\}$ define a probability distribution on $\{1,\dots,m\}$ (Dirac notation is implemented~\cite{Di39}). Set 
\[
\pi:=\sum_{i,j}\sqrt{p_{i}p_{j}}|i\>\<j|\otimes|i\>\<j|\in\Alg{A}\otimes\Alg{A}.
\]
Then $\pi$ is a purification of $\alpha$ (cf.\ \cite[Section 2.5]{NiCh11}). 
Notice also that if $P_{\alpha}$ is the \define{support projection} of $\alpha$, uniquely defined as the minimal orthogonal projection $P_{\alpha}\in\matr_{m}$ (meaning $P_{\alpha}^{\dag}P_{\alpha}=P_{\alpha}$) such that $P_{\alpha}\alpha=\alpha=\alpha P_{\alpha}$, then 
\be
\label{eq:supportprojpurification}
(\Ad_{P_{\alpha}}\otimes\id_{\matr_{m}})\circ\pi=\pi=(\id_{\matr_{m}}\otimes\Ad_{P_{\alpha}})\circ\pi, 
\ee
where $\Ad_{V}:\Alg{A}\to\Alg{A}$ is the CP map sending $A\in\Alg{A}$ to $VAV^{\dag}$. 

\item
Again, say $\Alg{A}=\matr_{m}$, and let $\Alg{A}\xrightarrow{\map{G}}\Alg{E}$ be any morphism in $\NCFinStoch$, i.e., a CPTP map. Write $\alpha=\sum_{i=1}^{m}p_{i}|i\>\<i|$ as in part~(\ref{item:purification}). Set $\pi_{\map{G}}:=\sum_{i,j}\sqrt{p_{i}p_{j}}|i\>\<j|\otimes \map{G}\big(|i\>\<j|\big)\in\Alg{A}\otimes\Alg{E}$. Then $\pi_{\map{G}}$ is a dilation of $\alpha$. Note that $\pi_{\map{G}}$ is not, in general, a purification of $\alpha$. 

\item
\label{item:purificationA}
Now let $\Alg{A}$ be an arbitrary finite-dimensional $C^*$-algebra. Then there exists a Hilbert space $\Hi$ of dimension $\sqrt{\Tr[1_{\Alg{A}}]}$ together with an injective trace-preserving unital $*$-homomorphism $\Alg{A}\xrightarrow{\varphi}\mathcal{L}(\Hi)$ (which is, in particular, CPTP) and a CPTP unital surjection $\mathcal{L}(\Hi)\xrightarrow{\psi}\Alg{A}$ such that $\psi\circ\varphi=\id_{\Alg{A}}$~\cite{Fa01}. This last equation says that $\psi$ is a \define{conditional expectation} of $\varphi$ if $\Alg{A}$ is viewed as a subalgebra of $\mathcal{L}(\Hi)$ via the map $\varphi$~\cite{Um54,GPRR21}. Given a state $\C\xrightarrow{\alpha}\Alg{A}$, we can push this state forward to a state $\C\xrightarrow{\varphi\circ\alpha}\mathcal{L}(\Hi)$. Write $\varphi(\alpha)=\sum_{i=1}^{m}p_{i}|i\>\<i|$ in terms of an orthonormal basis of eigenvectors in $\Hi$ of $\varphi(\alpha)$ so that the $p_{i}$ define a probability distribution on $\{1,\dots,m\}$. Set 
\[
\nu:=\sum_{i,j}\sqrt{p_i p_j}|i\>\<j|\otimes|i\>\<j|\in\mathcal{L}(\Hi)\otimes\mathcal{L}(\Hi).
\] 
Then $\nu$ is a purification of $\varphi\circ\alpha$ and 
\[
\pi:=(\psi\otimes\id_{\mathcal{L}(\Hi)})\circ\nu
\]
is a dilation of $\alpha$. Indeed, the dilation condition follows from 
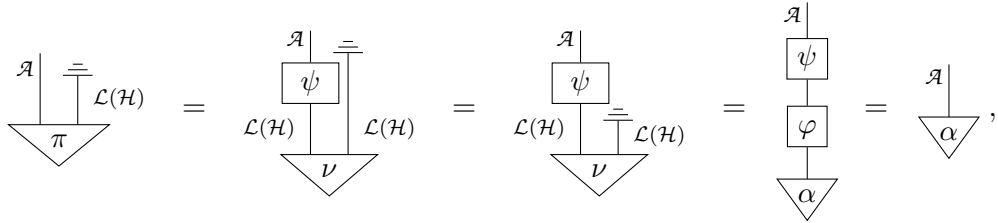
\[
\vcenter{\hbox{%
\begin{tikzpicture}[font=\small]
\node[state] (omega) at (0,0) {\;$\pi$\;};
\coordinate (X) at (-0.25,0.95) {};
\node[discarder] (Y) at (0.25,0.65) {};
\draw (omega) ++(-0.25, 0) to (X);
\draw (omega) ++(0.25, 0) to (Y);
\path[scriptstyle]
node at (-0.45,0.75) {$\Alg{A}$}
node at (0.8,0.35) {$\mathcal{L}(\Hi)$};
\end{tikzpicture}}}
\;\;=\;\;
\vcenter{\hbox{%
\begin{tikzpicture}[font=\small]
\node[state] (omega) at (0,0) {\;$\nu$\;};
\coordinate (X) at (-0.25,1.65) {};
\node[discarder] (Y) at (0.25,1.35) {};
\node[arrow box] (h) at (-0.25,0.95) {\;$\psi$\;};
\draw (omega) ++(-0.25, 0) to (h);
\draw (h) to (X);
\draw (omega) ++(0.25, 0) to (Y);
\path[scriptstyle]
node at (-0.45,1.55) {$\Alg{A}$}
node at (-0.8,0.35) {$\mathcal{L}(\Hi)$}
node at (0.8,0.35) {$\mathcal{L}(\Hi)$};
\end{tikzpicture}}}
\;\;=\;\;
\vcenter{\hbox{%
\begin{tikzpicture}[font=\small]
\node[state] (omega) at (0,0) {\;$\nu$\;};
\coordinate (X) at (-0.25,1.65) {};
\node[discarder] (Y) at (0.25,0.45) {};
\node[arrow box] (h) at (-0.25,0.95) {\;$\psi$\;};
\draw (omega) ++(-0.25, 0) to (h);
\draw (h) to (X);
\draw (omega) ++(0.25, 0) to (Y);
\path[scriptstyle]
node at (-0.45,1.55) {$\Alg{A}$}
node at (-0.8,0.35) {$\mathcal{L}(\Hi)$}
node at (0.8,0.25) {$\mathcal{L}(\Hi)$};
\end{tikzpicture}}}
\;\;=\;\;
\vcenter{\hbox{%
\begin{tikzpicture}[font=\small]
\node[state] (omega) at (0,-0.8) {$\alpha$};
\node[arrow box] (q) at (0,-0.1) {$\varphi$};
\node[arrow box] (g) at (0,0.8) {$\psi$};
\draw (omega) to (q);
\draw (q) to (g);
\draw (g) to (0,1.55);
\path[scriptstyle]
node at (-0.20,1.40) {$\Alg{A}$};
\end{tikzpicture}}}
\;\;=\;\;
\vcenter{\hbox{%
\begin{tikzpicture}[font=\small]
\node[state] (omega) at (0,-0.8) {$\alpha$};
\coordinate (q) at (0,-0.1);
\draw (omega) to (q);
\path[scriptstyle]
node at (-0.20,-0.25) {$\Alg{A}$};
\end{tikzpicture}}}
\;, 
\]
which, using standard algebraic notation rather than string diagrams, reads
\[
\Tr_{\mathcal{L}(\Hi)}\circ\pi
=\Tr_{\mathcal{L}(\Hi)}\circ\big((\psi\otimes\id_{\mathcal{L}(\Hi)})\circ\nu\big)
=\psi\circ\big(\Tr_{\mathcal{L}(\Hi)}\circ\nu\big)
=\psi\circ(\varphi\circ\alpha)
=(\psi\circ\varphi)\circ\alpha
=\alpha.
\]
Note that this reduces to the special case in Example~\ref{ex:dilationsCAlg} part~(\ref{item:purification}) when $\Alg{A}=\matr_{m}$ by taking $\Hi=\C^{m}$ and using the identity for the injective $*$-homomorphism $\varphi$. 
However, in the more general case when $\Alg{A}$ is not necessarily a matrix algebra, $\pi$ is not necessarily a pure state as an element of $\Alg{A}\otimes\mathcal{L}(\Hi)$. A simple example is given by $\Alg{A}=\C^{\{0,1\}}$ with the state $\alpha=\frac{1}{2}1_{\C^{\{0,1\}}}$. Although 
\[
\nu=\frac{1}{2}\Big(|0\>\<0|\otimes|0\>\<0|+|0\>\<1|\otimes|0\>\<1|+|1\>\<0|\otimes|1\>\<0|+|1\>\<1|\otimes|1\>\<1|\Big)
\]
is pure as an element of $\mathcal{L}(\C^2)\otimes\mathcal{L}(\C^2)$, applying $\psi\otimes\id_{\mathcal{L}(\C^2)}$ to it gives 
\[
\pi=\frac{1}{2}\delta_{0}\otimes |0\>\<0|+\frac{1}{2}\delta_{1}\otimes|1\>\<1|,
\]
where $\delta_{x}\in\C^{\{0,1\}}$ is the Kronecker delta function at $x\in\{0,1\}$. This shows that $\pi\in\C^{\{0,1\}}\otimes\mathcal{L}(\C^2)$ is not a pure state, since it is a non-trivial convex combination of two pure states. 

\item
Given the same setup as in the preceding example (Example~\ref{ex:dilationsCAlg} part~(\ref{item:purificationA})), it is also the case that $(\id_{\Alg{A}}\otimes\map{G})\circ\pi=(\psi\otimes\map{G})\circ\nu$ is a dilation of $\alpha$ for all CPTP maps $\mathcal{L}(\Hi)\xrightarrow{\map{G}}\Alg{E}$.
\end{enumerate}
\ex

\br
Although the GNS construction~\cite{GN43,Se47,PaGNS} (or more generally Stinespring's construction~\cite{St55,Pa18}) provides a representation of an arbitrary state $\alpha\in\Alg{A}$ in terms of a pure state in a Hilbert space $\mathcal{K}$, the algebra $\mathcal{L}(\mathcal{K})$ is not in general in the tensor product form $\Alg{A}\otimes\Alg{E}$ for some algebra $\Alg{E}$ and is therefore not a dilation in the sense defined here. The definition of purification here in terms of a dilation has been used as an axiom to isolate quantum theory among other GPTs in Ref.~\cite{CDP10}. 
\er

\blem
\label{lem:quantumdilationcharacterization}
Let $\alpha\in\Alg{A}$ be a state on a finite-dimensional $C^*$-algebra $\Alg{A}$ and let $\C\xrightarrow{\pi'}\Alg{A}\otimes\Alg{E}$ be any dilation of $\alpha$. Then there exists a CPTP map $\mathcal{L}(\Hi)\xrightarrow{\map{G}}\Alg{E}$ such that $(\id_{\Alg{A}}\otimes\map{G})\circ\pi=\pi'$, where $\Hi$ and $\pi$ are as in Example~\ref{ex:dilationsCAlg} part~(\ref{item:purificationA}). 
\elem

This is shown in~\cite[Theorem 2.2]{CHM22} in the case where $\Alg{A}$ is a matrix algebra. Because there are additional steps needed for the proof to extend to arbitrary algebras, we include a full proof here.

\bprf
Following the construction as in Example~\ref{ex:dilationsCAlg} part~(\ref{item:purificationA}), but for $\Alg{E}$, let $\mathcal{K}$ be a Hilbert space of dimension $\sqrt{\Tr[1_{\Alg{E}}]}$, let $\Alg{E}\xrightarrow{\eta}\mathcal{L}(\mathcal{K})$ be an injective trace-preserving unital $*$-homomorphism with associated conditional expectation $\mathcal{L}(\mathcal{K})\xrightarrow{\zeta}\Alg{E}$. Let 
\[
\C\xrightarrow{\nu'}\mathcal{L}(\Hi\otimes\mathcal{K})\otimes\mathcal{L}(\Hi\otimes\mathcal{K})\cong\mathcal{L}(\Hi)\otimes\mathcal{L}(\mathcal{K})\otimes\mathcal{L}(\Hi)\otimes\mathcal{L}(\mathcal{K})
\]
be a purification of $(\varphi\otimes\eta)\circ\pi'$, constructed analogously to Example~\ref{ex:dilationsCAlg} part~\ref{item:purification}. Then $\nu'$ is also a purification of $\varphi\circ\alpha$ because $\nu'$ is pure by construction and because
\[
\vcenter{\hbox{%
\begin{tikzpicture}[font=\small]
\node[state] (omega) at (0,0) {\;\;$\nu'$\;\;};
\coordinate (1) at (-0.75,0.95) {};
\node[discarder] (2) at (-0.25,0.65) {};
\node[discarder] (3) at (0.25,0.65) {};
\node[discarder] (4) at (0.75,0.65) {};
\draw (omega) ++(-0.75, 0) to (1);
\draw (omega) ++(-0.25, 0) to (2);
\draw (omega) ++(0.25, 0) to (3);
\draw (omega) ++(0.75, 0) to (4);
\path[scriptstyle]
node at (-1.25,0.75) {$\mathcal{L}(\Hi)$};
\end{tikzpicture}}}
\;\;=\;\;
\vcenter{\hbox{%
\begin{tikzpicture}[font=\small]
\node[state] (omega) at (0,0) {\;\;$\pi'$\;\;};
\coordinate (X) at (-0.5,1.65) {};
\node[discarder] (Y) at (0.5,1.35) {};
\node[arrow box] (h) at (-0.5,0.85) {\;$\varphi$\;};
\node[arrow box] (g) at (0.5,0.85) {\;$\eta$\;};
\draw (omega) ++(-0.5, 0) to (h);
\draw (h) to (X);
\draw (omega) ++(0.5, 0) to (g);
\draw (g) to (Y);
\path[scriptstyle]
node at (-1.00,1.45) {$\mathcal{L}(\Hi)$}
node at (-0.7,0.35) {$\Alg{A}$}
node at (0.7,0.35) {$\Alg{E}$};
\end{tikzpicture}}}
\;\;=\;\;
\vcenter{\hbox{%
\begin{tikzpicture}[font=\small]
\node[state] (omega) at (0,0) {\;$\pi'$\;};
\coordinate (X) at (-0.25,1.65) {};
\node[discarder] (Y) at (0.25,0.35) {};
\node[arrow box] (h) at (-0.25,0.85) {\;$\varphi$\;};
\draw (omega) ++(-0.25, 0) to (h);
\draw (h) to (X);
\draw (omega) ++(0.25, 0) to (Y);
\path[scriptstyle]
node at (-0.75,1.45) {$\mathcal{L}(\Hi)$}
node at (-0.45,0.35) {$\Alg{A}$};
\end{tikzpicture}}}
\;\;=\;\;
\vcenter{\hbox{%
\begin{tikzpicture}[font=\small]
\node[state] (omega) at (0,-0.8) {$\alpha$};
\node[arrow box] (q) at (0,-0.1) {$\varphi$};
\draw (omega) to (q);
\draw (q) to (0,0.65);
\path[scriptstyle]
node at (-0.45,0.45) {$\mathcal{L}(\Hi)$};
\end{tikzpicture}}}
\;,
\]
i.e., 
\[
\begin{split}
\big(\id_{\mathcal{L}(\Hi)}\otimes\Tr_{\mathcal{L}(\mathcal{K})}\otimes\Tr_{\mathcal{L}(\Hi\otimes\mathcal{K})}\big)\circ\nu'
&=\big(\id_{\mathcal{L}(\Hi)}\otimes\Tr_{\mathcal{L}(\mathcal{K})}\big)\circ\big((\varphi\otimes\eta)\circ\pi'\big)\\
&=\varphi\circ(\Tr_{\mathcal{L}(\mathcal{K})}\circ\pi')
=\varphi\circ\alpha,
\end{split}
\]
which shows that the dilation condition holds. Thus, by the essential uniqueness of purifications (see \cite[Section 2 in Preliminaries]{Ho21} and \cite[Theorem~5.1]{Wilde2013} for example), there exists an isometry $W:\Hi\to\mathcal{K}\otimes\Hi\otimes\mathcal{K}$ such that 
\[
(\id_{\mathcal{L}(\Hi)}\otimes\Ad_{W})\circ\nu=\nu',
\]
i.e., 
\[
\vcenter{\hbox{%
\begin{tikzpicture}[font=\small,xscale=-1]
\node[state] (omega) at (0,0) {\;$\nu$\;};
\coordinate (X) at (-0.35,1.65) {};
\coordinate (Y) at (0.35,1.65) {};
\node[arrow box] (h) at (-0.35,0.85) {$\Ad_{W}$};
\draw (omega) ++(-0.35, 0) to (h);
\draw (h) to (X);
\draw (omega) ++(0.35, 0) to (Y);
\path[scriptstyle]
node at (-1.45,1.45) {$\mathcal{L}(\mathcal{K}\otimes\Hi\otimes\mathcal{K})$}
node at (0.85,1.45) {$\mathcal{L}(\Hi)$}
node at (-0.85,0.35) {$\mathcal{L}(\Hi)$};
\end{tikzpicture}}}
\;\;=\;\;
\vcenter{\hbox{%
\begin{tikzpicture}[font=\small]
\node[state] (omega) at (0,0) {\;\;$\nu'$\;\;};
\coordinate (1) at (-0.75,0.95) {};
\coordinate (2) at (-0.25,0.95) {};
\coordinate (3) at (0.25,0.95) {};
\coordinate (4) at (0.75,0.95) {};
\draw (omega) ++(-0.75, 0) to (1);
\draw (omega) ++(-0.25, 0) to (2);
\draw (omega) ++(0.25, 0) to (3);
\draw (omega) ++(0.75, 0) to (4);
\end{tikzpicture}}}
\;.
\]
Let 
\[
\map{G}:=\big(\zeta\otimes\Tr_{\mathcal{L}(\Hi\otimes\mathcal{K})}\big)\circ\Ad_{W},
\]
i.e., 
\[
\vcenter{\hbox{%
\begin{tikzpicture}[font=\small]
\node[arrow box] (c) at (0,0) {$\map{G}$};
\draw (c) to (0,0.9);
\draw (c) to (0,-0.9);
\node at (-0.5,-0.70) {\scriptsize $\mathcal{L}(\Hi)$};
\node at (-0.2,0.70) {\scriptsize $\Alg{E}$};
\end{tikzpicture}}}
\;\;:=\;\;
\vcenter{\hbox{%
\begin{tikzpicture}[font=\small]
\coordinate (X) at (0,-1.0) {};
\node[arrow box] (g) at (-0.7,1.0) {$\zeta$};
\node[discarder] (dA) at (0,1.0) {};
\node[discarder] (dE) at (0.7,1.0) {};
\node[arrow box] (AdW) at (0,0) {\;\;\;\;\;$\Ad_{W}$\;\;\;\;\;};
\coordinate (Y) at (-0.7,1.8);
\draw (X) to (AdW);
\draw (AdW) ++(-0.7, 0.25) to (g);
\draw (AdW) ++(0, 0.25) to (dA);
\draw (AdW) ++(0.7, 0.25) to (dE);
\draw (g) to (Y);
\path[scriptstyle]
node at (-0.5,-0.75) {$\mathcal{L}(\Hi)$}
node at (-0.95,1.65) {$\Alg{E}$};
\end{tikzpicture}}}
\;,
\]
which is a morphism in $\NCFinStoch$ because $\Ad_{W}$ is CPTP by the fact that $W$ is an isometry. Therefore, 
\[
\vcenter{\hbox{%
\begin{tikzpicture}[font=\small]
\node[state] (omega) at (0,0) {\;$\pi$\;};
\coordinate (X) at (0.25,1.65) {};
\coordinate (Y) at (-0.25,1.65) {};
\node[arrow box] (h) at (0.25,0.85) {\;$\map{G}$\;};
\draw (omega) ++(0.25, 0) to (h);
\draw (h) to (X);
\draw (omega) ++(-0.25, 0) to (Y);
\path[scriptstyle]
node at (0.45,1.45) {$\Alg{E}$}
node at (-0.45,1.45) {$\Alg{A}$}
node at (0.75,0.3) {$\mathcal{L}(\Hi)$};
\end{tikzpicture}}}
\;\;=\;\;
\vcenter{\hbox{%
\begin{tikzpicture}[font=\small]
\node[state] (X) at (-0.3,-1.0) {\;\;\;\;$\nu$\;\;\;\;};
\node[arrow box] (g) at (-0.7,1.0) {$\zeta$};
\node[arrow box] (h) at (-1.4,1.0) {$\psi$};
\node[discarder] (dA) at (0,1.0) {};
\node[discarder] (dE) at (0.7,1.0) {};
\node[arrow box] (AdW) at (0,0) {\;\;\;\;\;$\Ad_{W}$\;\;\;\;\;};
\coordinate (Y) at (-0.7,1.8);
\coordinate (Z) at (-1.4,1.8);
\draw (X) ++(0.3, 0) to (AdW);
\draw (X) ++(-1.1,0) to (h);
\draw (AdW) ++(-0.7, 0.25) to (g);
\draw (AdW) ++(0, 0.25) to (dA);
\draw (AdW) ++(0.7, 0.25) to (dE);
\draw (g) to (Y);
\draw (h) to (Z);
\path[scriptstyle]
node at (-0.9,-0.65) {$\mathcal{L}(\Hi)$}
node at (0.5,-0.65) {$\mathcal{L}(\Hi)$}
node at (-0.55,1.65) {$\Alg{E}$}
node at (-1.65,1.65) {$\Alg{A}$};
\end{tikzpicture}}}
\;\;=\;\;
\vcenter{\hbox{%
\begin{tikzpicture}[font=\small]
\node[state] (omega) at (0,0) {\;\;\;$\nu'$\;\;\;};
\node[arrow box] (g) at (-0.95,1.0) {$\psi$};
\node[arrow box] (h) at (-0.3,1.0) {$\zeta$};
\coordinate (1) at (-0.95,1.95) {};
\coordinate (2) at (-0.3,1.95) {};
\node[discarder] (3) at (0.3,0.95) {};
\node[discarder] (4) at (0.95,0.95) {};
\draw (omega) ++(-0.95, 0) to (g);
\draw (omega) ++(-0.3, 0) to (h);
\draw (g) to (1);
\draw (h) to (2);
\draw (omega) ++(0.3, 0) to (3);
\draw (omega) ++(0.95, 0) to (4);
\path[scriptstyle]
node at (-1.15,1.8) {$\Alg{A}$}
node at (-0.15,1.8) {$\Alg{E}$};
\end{tikzpicture}}}
\;\;=\;\;
\vcenter{\hbox{%
\begin{tikzpicture}[font=\small]
\node[state] (omega) at (0,0) {\;\;$\pi'$\;\;};
\coordinate (X) at (-0.5,2.55) {};
\coordinate (Y) at (0.5,2.55) {};
\node[arrow box] (g2) at (0.5,1.85) {\;$\zeta$\;};
\node[arrow box] (h) at (-0.5,0.85) {\;$\varphi$\;};
\node[arrow box] (h2) at (-0.5,1.85) {\;$\psi$\;};
\node[arrow box] (g) at (0.5,0.85) {\;$\eta$\;};
\draw (omega) ++(-0.5, 0) to (h);
\draw (h) to (h2);
\draw (h2) to (X);
\draw (omega) ++(0.5, 0) to (g);
\draw (g) to (g2);
\draw (g2) to (Y);
\path[scriptstyle]
node at (-1.00,1.35) {$\mathcal{L}(\Hi)$}
node at (1.00,1.35) {$\mathcal{L}(\mathcal{K})$}
node at (-0.7,0.35) {$\Alg{A}$}
node at (0.7,0.35) {$\Alg{E}$}
node at (-0.7,2.45) {$\Alg{A}$}
node at (0.7,2.45) {$\Alg{E}$};
\end{tikzpicture}}}
\;\;=\;\;
\vcenter{\hbox{%
\begin{tikzpicture}[font=\small]
\node[state] (omega) at (0,0) {\;$\pi'$\;};
\coordinate (X) at (-0.25,0.95) {};
\coordinate (Y) at (0.25,0.95) {};
\draw (omega) ++(-0.25, 0) to (X);
\draw (omega) ++(0.25, 0) to (Y);
\path[scriptstyle]
node at (-0.45,0.75) {$\Alg{A}$}
node at (0.45,0.75) {$\Alg{E}$};
\end{tikzpicture}}}
\;,
\]
i.e., 
\[
\begin{split}
(\id_{\Alg{A}}\otimes\map{G})\circ\pi
&=(\psi\otimes\zeta\otimes\Tr_{\mathcal{L}(\Hi\otimes\mathcal{K})})\circ(\id_{\mathcal{L}(\Hi)}\otimes\Ad_{W})\circ\nu
=(\psi\otimes\zeta\otimes\Tr_{\mathcal{L}(\Hi\otimes\mathcal{K})})\circ\nu'\\
&=(\psi\otimes\zeta)\circ(\varphi\otimes\eta)\circ\pi'
=\big((\psi\circ\varphi)\otimes(\zeta\circ\eta)\big)\circ\pi'
=\pi',
\end{split}
\]
which proves the claim. 
\eprf

\bd
\label{defn:dilationalequality}
Let $\mathbf{C}$ be a semicartesian category, let $\map{P}:\Alg{Z}\to\Alg{A}$ and $\map{E},\map{F}:\Alg{A}\to\Alg{B}$ be morphisms in $\mathbf{C}$. Then $\map{E}$ is said to be \define{$\map{P}$-dilationally equal to} $\map{F}$, written%
\footnote{In this paper, the notation $\aeequals{\map{P}}$ will be used exclusively to refer to dilational equality, which differs in general from a.e.\ equivalence as defined in Refs.~\cite{ChJa18,PaBayes} (see Remark~\ref{rmk:aeequalsvsdilationalequality} for more details).}
 $\map{E}\aeequals{\map{P}}\map{F}$, iff $(\map{E}\otimes\id_{\Alg{E}})\circ\map{D}=(\map{F}\otimes\id_{\Alg{E}})\circ\map{D}$ for every dilation $\map{D}:\Alg{Z}\to\Alg{A}\otimes\Alg{E}$ of $\map{P}$, i.e., 
\[
\vcenter{\hbox{%
\begin{tikzpicture}[font=\small]
\node[arrow box] (c) at (0,0) {$\map{E}$};
\draw (c) to (0,0.8);
\draw (c) to (0,-0.8);
\node at (-0.2,-0.70) {\scriptsize $\Alg{A}$};
\node at (-0.2,0.70) {\scriptsize $\Alg{B}$};
\end{tikzpicture}}}
\;\;\;\aeequals{\map{P}}\;\;
\vcenter{\hbox{%
\begin{tikzpicture}[font=\small]
\node[arrow box] (c) at (0,0) {$\map{F}$};
\draw (c) to (0,0.8);
\draw (c) to (0,-0.8);
\node at (-0.2,-0.70) {\scriptsize $\Alg{A}$};
\node at (-0.2,0.70) {\scriptsize $\Alg{B}$};
\end{tikzpicture}}}
\;\;\;\iff\;\;\;
\vcenter{\hbox{%
\begin{tikzpicture}[font=\small]
\node[arrow box] (D) at (0,-0.1) {\;\;\;$\map{D}$\;\;\;};
\coordinate (O) at (0,-0.8);
\coordinate (X) at (-0.35,1.65) {};
\coordinate (Y) at (0.35,1.65) {};
\node[arrow box] (h) at (-0.35,0.85) {$\map{E}$};
\draw (O) to (D);
\draw (D) ++(-0.35, 0.25) to (h);
\draw (h) to (X);
\draw (D) ++(0.35, 0.25) to (Y);
\path[scriptstyle]
node at (-0.2,-0.65) {$\Alg{Z}$}
node at (-0.55,1.45) {$\Alg{B}$}
node at (0.55,1.45) {$\Alg{E}$}
node at (-0.55,0.4) {$\Alg{A}$};
\end{tikzpicture}}}
\;\;=\;\;
\vcenter{\hbox{%
\begin{tikzpicture}[font=\small]
\node[arrow box] (D) at (0,-0.1) {\;\;\;$\map{D}$\;\;\;};
\coordinate (O) at (0,-0.8);
\coordinate (X) at (-0.35,1.65) {};
\coordinate (Y) at (0.35,1.65) {};
\node[arrow box] (h) at (-0.35,0.85) {$\map{F}$};
\draw (O) to (D);
\draw (D) ++(-0.35, 0.25) to (h);
\draw (h) to (X);
\draw (D) ++(0.35, 0.25) to (Y);
\path[scriptstyle]
node at (-0.2,-0.65) {$\Alg{Z}$}
node at (-0.55,1.45) {$\Alg{B}$}
node at (0.55,1.45) {$\Alg{E}$}
node at (-0.55,0.4) {$\Alg{A}$};
\end{tikzpicture}}}
\quad\text{for all 
$
\vcenter{\hbox{
\begin{tikzpicture}[font=\small]
\node[arrow box] (p) at (0,0) {\;\;$\map{D}$\;\;};
\coordinate (X) at (-0.25,0.8);
\coordinate (Y) at (0.25,0.8);
\draw (0,-0.8) to (p);
\draw (p.north)++(-0.25,0) to (X);
\draw (p.north)++(0.25,0) to (Y);
\node at (-0.2,-0.70) {\scriptsize $\Alg{Z}$};
\node at (-0.4,0.70) {\scriptsize $\Alg{A}$};
\node at (0.4,0.70) {\scriptsize $\Alg{E}$};
\end{tikzpicture}}}
$
such that 
$
\vcenter{\hbox{
\begin{tikzpicture}[font=\small]
\node[arrow box] (p) at (0,0) {\;\;$\map{D}$\;\;};
\coordinate (X) at (-0.25,0.8);
\node[discarder] (Y) at (0.25,0.5) {};
\draw (0,-0.8) to (p);
\draw (p.north)++(-0.25,0) to (X);
\draw (p.north)++(0.25,0) to (Y);
\node at (-0.2,-0.70) {\scriptsize $\Alg{Z}$};
\node at (-0.4,0.70) {\scriptsize $\Alg{A}$};
\end{tikzpicture}}}
\;\;=\;\;
\vcenter{\hbox{%
\begin{tikzpicture}[font=\small]
\node[arrow box] (c) at (0,0) {$\map{P}$};
\draw (c) to (0,0.8);
\draw (c) to (0,-0.8);
\node at (-0.2,-0.70) {\scriptsize $\Alg{Z}$};
\node at (-0.2,0.70) {\scriptsize $\Alg{A}$};
\end{tikzpicture}}}
$
\;.
}
\]
\ed

\bx
In the semicartesian category $\FinStoch$, dilational equality coincides with a.e.\ equality as defined by Ref.~\cite{ChJa18}, as shown more generally for arbitrary causal Markov categories in~\cite[Proposition 4.2]{FGGPS22}. 
\ex

In the quantum setting, we have the following result%
\footnote{I thank Tobias Fritz and Tom{\'a}{\v s} Gonda for helpful discussions related to this.}.

\bn
\label{prop:dilationalequalityQuantum}
Let $\C\xrightarrow{\alpha}\Alg{A}$ be a state in $\NCFinStoch$. Then two morphisms $\map{E},\map{F}:\Alg{A}\to\Alg{B}$ are $\alpha$-dilationally equal if and only if $\map{E}\circ\Ad_{P_{\alpha}}=\map{F}\circ\Ad_{P_{\alpha}}$, where $P_{\alpha}$ is the support projection of $\alpha$, uniquely defined as the minimal orthogonal projection $P_{\alpha}\in\Alg{A}$ (meaning $P_{\alpha}^{\dag}P_{\alpha}=P_{\alpha}$) such that $P_{\alpha}\alpha=\alpha=\alpha P_{\alpha}$. 
\en

\bprf
First assume that $\map{E},\map{F}:\Alg{A}\to\Alg{B}$ are $\alpha$-dilationally equal. Let $\Hi, \varphi, \psi, \nu$, and $\pi$ be as in Example~\ref{ex:dilationsCAlg} part~(\ref{item:purificationA}). Then, since $\pi$ is a dilation of $\alpha$, 
\[
\Big((\map{E}\circ\psi)\otimes\id_{\mathcal{L}(\Hi)}\Big)\left(\sum_{i,j}\sqrt{p_{i}p_{j}}|i\>\<j|\otimes|i\>\<j|\right)
=\Big((\map{F}\circ\psi)\otimes\id_{\mathcal{L}(\Hi)}\Big)\left(\sum_{i,j}\sqrt{p_{i}p_{j}}|i\>\<j|\otimes|i\>\<j|\right)
\]
by the assumption of $\alpha$-dilational equality. 
By linear independence of the $|i\>\<j|$ as elements of $\mathcal{L}(\Hi)$, it follows that each $i,j$ term in the summation must be equal, i.e., 
\be
\label{eq:aeEF}
\sqrt{p_{i}p_{j}}(\map{E}\circ\psi)\big(|i\>\<j|\big)
=\sqrt{p_{i}p_{j}}(\map{F}\circ\psi)\big(|i\>\<j|\big)
\ee
for all $i,j$. Let $P_{\alpha}\in\Alg{A}$ and $P_{\varphi(\alpha)}\in\mathcal{L}(\Hi)$ denote the support projections of $\alpha$ and $\varphi\circ\alpha$, respectively. By Equation~\eqref{eq:aeEF}, we obtain
\[
(\map{E}\circ\psi)\big(|i\>\<j|\big)=(\map{F}\circ\psi)\big(|i\>\<j|\big)
\]
for all $i,j$ such that $p_{i}\ne0$ and $p_{j}\ne0$, which is equivalent to 
\[
\map{E}\circ\psi\circ\Ad_{P_{\varphi(\alpha)}}=\map{F}\circ\psi\circ\Ad_{P_{\varphi(\alpha)}}. 
\]
Pre-composing both sides with $\varphi$ gives
\[
\xy0;/r.22pc/:
(-20,13)*+{\map{E}\circ\psi\circ\Ad_{P_{\varphi(\alpha)}}\circ\varphi}="1";
(-30,0)*+{\map{E}\circ\psi\circ\varphi\circ\Ad_{P_{\alpha}}}="2";
(-20,-13)*+{\map{E}\circ\Ad_{P_{\alpha}}}="3";
(20,-13)*+{\map{F}\circ\Ad_{P_{\alpha}}}="4";
(30,0)*+{\map{F}\circ\psi\circ\varphi\circ\Ad_{P_{\alpha}}}="5";
(20,13)*+{\map{F}\circ\psi\circ\Ad_{P_{\varphi(\alpha)}}\circ\varphi}="6";
{\ar@{=}@/_0.85pc/"1";"2"_{\text{}}};
{\ar@{=}@/_0.55pc/"2";"3"_(0.45){\text{}}};
{\ar@{=}@/_0.25pc/"6";"1"_(0.5){\text{}}};
{\ar@{=}@/_0.55pc/"4";"5"_(0.65){\text{}}};
{\ar@{=}@/_0.85pc/"5";"6"_{\text{}}};
\endxy
,
\]
where the equalities come from the properties of the support projections and the maps $\psi$ and $\varphi$. This proves one direction of the claim.  

Conversely, suppose $\map{E}\circ\Ad_{P_{\alpha}}=\map{F}\circ\Ad_{P_{\alpha}}$. Let $\pi':\C\to\Alg{A}\otimes\Alg{E}$ be a dilation of $\alpha$. By Lemma~\ref{lem:quantumdilationcharacterization}, there exists a CPTP map $\mathcal{L}(\Hi)\xrightarrow{\map{G}}\Alg{E}$ such that $(\id_{\Alg{A}}\otimes\map{G})\circ\pi=\pi'$, where $\Hi,\psi,\varphi,\pi$ define a purification of $\alpha$ as in Example~\ref{ex:dilationsCAlg} part~(\ref{item:purificationA}). Therefore, 
\begin{align*}
(\map{E}\otimes\id_{\Alg{E}})\circ\pi'&=(\map{E}\otimes\map{G})\circ\pi &&\mbox{by Lemma~\ref{lem:quantumdilationcharacterization}}\\
&=(\map{E}\otimes\map{G})\circ(\psi\otimes\id_{\mathcal{L}(\Hi)})\circ\nu&&\mbox{by definition of $\pi$}\\
&=(\map{E}\otimes\map{G})\circ\big((\psi\circ\Ad_{P_{\varphi(\alpha)}})\otimes\id_{\mathcal{L}(\Hi)}\big)\circ\nu&&\mbox{by the property~\eqref{eq:supportprojpurification}} \\
&=\big((\map{E}\circ\Ad_{P_{\alpha}}\circ\psi)\otimes\map{G}\big)\circ\nu&&\mbox{by conditional expectation property}\\
&=\big((\map{F}\circ\Ad_{P_{\alpha}}\circ\psi)\otimes\map{G}\big)\circ\nu&&\mbox{by assumption}\\
&=\cdots=(\map{F}\otimes\id_{\Alg{E}})\circ\pi'&&\mbox{by an analogous calculation in reverse},
\end{align*}
which, using string diagrams, reads
\begingroup
\allowdisplaybreaks
\begin{align*}
\vcenter{\hbox{%
\begin{tikzpicture}[font=\small]
\node[state] (omega) at (0,0) {\;$\pi'$\;};
\coordinate (X) at (-0.35,1.65) {};
\coordinate (Y) at (0.35,1.65) {};
\node[arrow box] (h) at (-0.35,0.85) {$\map{E}$};
\draw (omega) ++(-0.35, 0) to (h);
\draw (h) to (X);
\draw (omega) ++(0.35, 0) to (Y);
\path[scriptstyle]
node at (-0.55,1.45) {$\Alg{B}$}
node at (0.55,1.00) {$\Alg{E}$}
node at (-0.55,0.35) {$\Alg{A}$};
\end{tikzpicture}}}
\;\;&=\;\;
\vcenter{\hbox{%
\begin{tikzpicture}[font=\small]
\node[state] (omega) at (0,0) {\;\;$\pi$\;\;};
\coordinate (X) at (0.45,1.65) {};
\coordinate (Y) at (-0.45,1.65) {};
\node[arrow box] (e) at (-0.45,0.85) {\;$\map{E}$\;};
\node[arrow box] (g) at (0.45,0.85) {\;$\map{G}$\;};
\draw (omega) ++(0.45, 0) to (g);
\draw (g) to (X);
\draw (omega) ++(-0.45, 0) to (e);
\draw (e) to (Y);
\path[scriptstyle]
node at (0.65,1.45) {$\Alg{E}$}
node at (-0.65,0.35) {$\Alg{A}$}
node at (-0.65,1.45) {$\Alg{B}$}
node at (0.95,0.3) {$\mathcal{L}(\Hi)$};
\end{tikzpicture}}}
\;\;=\;\;
\vcenter{\hbox{%
\begin{tikzpicture}[font=\small]
\node[state] (omega) at (0,0) {\;\;$\nu$\;\;};
\coordinate (X) at (0.45,2.65) {};
\coordinate (Y) at (-0.45,2.65) {};
\node[arrow box] (psi) at (-0.45,0.85) {\;$\psi$\;};
\node[arrow box] (e) at (-0.45,1.85) {\;$\map{E}$\;};
\node[arrow box] (g) at (0.45,1.85) {\;$\map{G}$\;};
\draw (omega) ++(0.45, 0) to (g);
\draw (g) to (X);
\draw (omega) ++(-0.45, 0) to (psi);
\draw (psi) to (e);
\draw (e) to (Y);
\path[scriptstyle]
node at (0.65,2.45) {$\Alg{E}$}
node at (-0.65,1.35) {$\Alg{A}$}
node at (-0.65,2.45) {$\Alg{B}$}
node at (-0.95,0.3) {$\mathcal{L}(\Hi)$}
node at (0.95,0.3) {$\mathcal{L}(\Hi)$};
\end{tikzpicture}}}
\;\;=\;\;
\vcenter{\hbox{%
\begin{tikzpicture}[font=\small]
\node[state] (omega) at (0,0) {\;\;$\nu$\;\;};
\coordinate (X) at (0.45,3.65) {};
\coordinate (Y) at (-0.45,3.65) {};
\node[arrow box] (W) at (-0.45,0.85) {$\Ad_{P_{\varphi(\alpha)}}$};
\node[arrow box] (psi) at (-0.45,1.85) {\;$\psi$\;};
\node[arrow box] (e) at (-0.45,2.85) {\;$\map{E}$\;};
\node[arrow box] (g) at (0.45,2.85) {\;$\map{G}$\;};
\draw (omega) ++(0.45, 0) to (g);
\draw (g) to (X);
\draw (omega) ++(-0.45, 0) to (W);
\draw (W) to (psi);
\draw (psi) to (e);
\draw (e) to (Y);
\path[scriptstyle]
node at (0.65,3.45) {$\Alg{E}$}
node at (-0.65,2.35) {$\Alg{A}$}
node at (-0.65,3.45) {$\Alg{B}$}
node at (-0.95,1.35) {$\mathcal{L}(\Hi)$}
node at (-0.95,0.3) {$\mathcal{L}(\Hi)$}
node at (0.95,0.3) {$\mathcal{L}(\Hi)$};
\end{tikzpicture}}}
\;\;=\;\;
\vcenter{\hbox{%
\begin{tikzpicture}[font=\small]
\node[state] (omega) at (0,0) {\;\;$\nu$\;\;};
\coordinate (X) at (0.45,3.65) {};
\coordinate (Y) at (-0.45,3.65) {};
\node[arrow box] (psi) at (-0.45,1.85) {$\Ad_{P_{\alpha}}$};
\node[arrow box] (W) at (-0.45,0.85) {\;$\psi$\;};
\node[arrow box] (e) at (-0.45,2.85) {\;$\map{E}$\;};
\node[arrow box] (g) at (0.45,2.85) {\;$\map{G}$\;};
\draw (omega) ++(0.45, 0) to (g);
\draw (g) to (X);
\draw (omega) ++(-0.45, 0) to (W);
\draw (W) to (psi);
\draw (psi) to (e);
\draw (e) to (Y);
\path[scriptstyle]
node at (0.65,3.45) {$\Alg{E}$}
node at (-0.65,2.35) {$\Alg{A}$}
node at (-0.65,3.45) {$\Alg{B}$}
node at (-0.65,1.35) {$\Alg{A}$}
node at (-0.95,0.3) {$\mathcal{L}(\Hi)$}
node at (0.95,0.3) {$\mathcal{L}(\Hi)$};
\end{tikzpicture}}}
\\
\;\;&=\;\;
\vcenter{\hbox{%
\begin{tikzpicture}[font=\small]
\node[state] (omega) at (0,0) {\;\;$\nu$\;\;};
\coordinate (X) at (0.45,3.65) {};
\coordinate (Y) at (-0.45,3.65) {};
\node[arrow box] (psi) at (-0.45,1.85) {$\Ad_{P_{\alpha}}$};
\node[arrow box] (W) at (-0.45,0.85) {\;$\psi$\;};
\node[arrow box] (e) at (-0.45,2.85) {\;$\map{F}$\;};
\node[arrow box] (g) at (0.45,2.85) {\;$\map{G}$\;};
\draw (omega) ++(0.45, 0) to (g);
\draw (g) to (X);
\draw (omega) ++(-0.45, 0) to (W);
\draw (W) to (psi);
\draw (psi) to (e);
\draw (e) to (Y);
\path[scriptstyle]
node at (0.65,3.45) {$\Alg{E}$}
node at (-0.65,2.35) {$\Alg{A}$}
node at (-0.65,3.45) {$\Alg{B}$}
node at (-0.65,1.35) {$\Alg{A}$}
node at (-0.95,0.3) {$\mathcal{L}(\Hi)$}
node at (0.95,0.3) {$\mathcal{L}(\Hi)$};
\end{tikzpicture}}}
\;\;=\;\;\cdots\;\;=\;\;
\vcenter{\hbox{%
\begin{tikzpicture}[font=\small]
\node[state] (omega) at (0,0) {\;$\pi'$\;};
\coordinate (X) at (-0.35,1.65) {};
\coordinate (Y) at (0.35,1.65) {};
\node[arrow box] (h) at (-0.35,0.85) {$\map{F}$};
\draw (omega) ++(-0.35, 0) to (h);
\draw (h) to (X);
\draw (omega) ++(0.35, 0) to (Y);
\path[scriptstyle]
node at (-0.55,1.45) {$\Alg{B}$}
node at (0.55,1.45) {$\Alg{E}$}
node at (-0.55,0.35) {$\Alg{A}$};
\end{tikzpicture}}}
\;.
\end{align*}
\endgroup
Hence, $\map{E}$ and $\map{F}$ are $\alpha$-dilationally equal, thus completing the proof. 
\eprf

\br
\label{rmk:aeequalsvsdilationalequality}
The notion of dilational equality from Proposition~\ref{prop:dilationalequalityQuantum}, which stems from Definition~\ref{defn:dilationalequality} as introduced in Ref.~\cite{FGGPS22}, is not equivalent to the definition of a.e.\ equivalence given for quantum Markov categories as in Ref.~\cite{PaBayes}, the latter of which is based on the original categorical definition of a.e.\ equivalence of Ref.~\cite{ChJa18}. Indeed, if two quantum channels are a.e.\ equivalent as in Ref.~\cite{PaBayes}, then they are dilationally equal as in Definition~\ref{defn:dilationalequality} due to Proposition~\ref{prop:dilationalequalityQuantum}. However, the converse need not be true, as discussed in Refs.~\cite{PaRuBayes,GPRR21} (based on Proposition~\ref{prop:dilationalequalityQuantum}).
\er

\bd
Let $\mathbf{C}$ be a semicartesian category. A state $I\xrightarrow{\alpha}\Alg{A}$ is said to be \define{faithful} iff for every pair $\map{E},\map{F}:\Alg{A}\to\Alg{B}$ of $\alpha$-dilationally equal morphisms implies $\map{E}=\map{F}$. 
\ed

A faithful state was defined in this way using dilations in Ref.~\cite{Chiribella14dilation}. Proposition~\ref{prop:dilationalequalityQuantum}  justifies this definition for arbitrary finite-dimensional $C^*$-algebras, as the following example shows. 

\bx
\label{ex:faithfulstateNC}
In $\NCFinStoch$, a state $\alpha\in\Alg{A}$ is faithful if and only if its support projection $P_{\alpha}$ equals the unit $1_{\Alg{A}}$. If $\Alg{A}=\matr_{m}$ so that $\alpha$ is a density matrix, this means that zero is not an eigenvalue of $\alpha$, i.e., $\alpha$ has rank $m$. More generally, as every finite-dimensional $C^*$-algebra is isomorphic to one of the form $\bigoplus_{x\in X}\matr_{m_{x}}$ for some finite set $X$ and natural numbers $m_{x}$ with $x\in X$~\cite{Fa01}, a state on it can be expressed as $\bigoplus_{x\in X}p_{x}\rho_{x}$, with $\{p_{x}\}$ a probability measure on $X$ and $\rho_{x}\in\matr_{m_{x}}$ a density matrix for all $x\in X$~\cite{PaRu19,PaRuBayes}. This state is then faithful if and only if $p_{x}>0$ and $\rho_{x}$ has rank $m_{x}$ for all $x\in X$. 
\ex

\bn
\label{prop:faithfultensor}
In the semicartesian categories $\FinStoch$ and $\NCFinStoch$, the tensor product of faithful states is again faithful. 
\en

\bprf
This follows from Example~\ref{ex:faithfulstateNC} and Proposition~\ref{prop:dilationalequalityQuantum} because if $\C\xrightarrow{\alpha}\Alg{A}$ and $\C\xrightarrow{\alpha'}\Alg{A}'$ are states, then $P_{\alpha\otimes\alpha'}=P_{\alpha}\otimes P_{\alpha'}$ and $\Ad_{P_{\alpha\otimes\alpha'}}=\Ad_{P_{\alpha}}\otimes\Ad_{P_{\alpha'}}$ for the support projections $P_{\alpha}$ and $P_{\alpha}$ of $\alpha$ and $\alpha'$, respectively. 
\eprf

\bq
What additional categorical axioms can be imposed on a semicartesian category $\mathbf{C}$ so that the tensor product of faithful states is faithful? 
\eq

\bd
\label{defn:fstatesSCC}
Let $\mathbf{C}$ be a semicartesian category and let $\states^{+}(\mathbf{C})$ be the \define{category of faithful states} and \define{state-preserving morphisms}. In detail, the objects of $\states^{+}(\mathbf{C})$ are pairs $(\Alg{A},\alpha)$, with $\Alg{I}\xrightarrow{\alpha}\Alg{A}$ a faithful state. A morphism from $(\Alg{A},\alpha)$ to $(\Alg{B},\beta)$ in $\states^{+}(\mathbf{C})$ is a morphism $\Alg{A}\xrightarrow{\map{E}}\Alg{B}$ in $\mathbf{C}$ such that $\map{E}\circ\alpha=\beta$. In this case, $\alpha$ is called the \define{prior} and $\beta$ the $\define{prediction}$.  
\ed

\br
Note that $\mathscr{S}^{+}(\mathbf{C})$ is a subcategory of the coslice/under category $I\downarrow\mathbf{C}$, whose objects are \emph{all} pairs $(\Alg{A},\alpha)$ with $\Alg{A}$ an object of $\mathbf{C}$ and $I\xrightarrow{\alpha}\Alg{A}$ a state. The reason we isolate a subcategory of \emph{faithful} states is because retrodiction, defined in the next section, will be typically only uniquely specified for faithful states, but not all states. An alternative construction of a category of states and state-preserving morphisms (that will be described shortly) is to use all states, not just faithful ones, but set the morphisms to be dilationally equal equivalence classes of state-preserving morphisms. 
\er

\bn
The category of faithful states $\states^{+}(\mathbf{C})$ is a symmetric monoidal (in fact, semicartesian) category for $\mathbf{C}=\FinStoch$ and $\mathbf{C}=\NCFinStoch$.
\en

\bprf
This follows from Proposition~\ref{prop:faithfultensor}.
\eprf

\blem
\label{lem:dilationcompandtensor}
Let $\mathbf{C}$ be a semicartesian category. 
Fix states $I\xrightarrow{\alpha}\Alg{A}$, $I\xrightarrow{\alpha'}\Alg{A}'$, and $I\xrightarrow{\beta}\Alg{B}$ in $\mathbf{C}$. 
Let $\map{E}_{1},\map{E}_{2}:\Alg{A}\to\Alg{B}$, $\map{F}_{1},\map{F}_{2}:\Alg{B}\to\Alg{C}$, and $\map{E}_{1}',\map{E}_{2}':\Alg{A}'\to\Alg{B}'$ be morphisms in $\mathbf{C}$ such that $\map{E}_{1}\circ\alpha=\beta=\map{E}_{2}\circ\alpha$, $\map{F}_{1}\circ\beta=\map{F}_{2}\circ\beta$, and $\map{E}_{1}'\circ\alpha'=\map{E}_{2}'\circ\alpha'$. 
Furthermore, suppose that $\map{E}_{1}\aeequals{\alpha}\map{E}_{2}$, $\map{F}_{1}\aeequals{\beta}\map{F}_{2}$, and $\map{E}_{1}'\aeequals{\alpha'}\map{E}_{2}'$. 
Then $\map{F}_{1}\circ\map{E}_{1}\aeequals{\alpha}\map{F}_{2}\circ\map{E}_{2}$ and $\map{E}_{1}\otimes\map{E}_{1}'\aeequals{\alpha\otimes\alpha'}\map{E}_{2}\otimes\map{E}_{2}'$.
\elem

\bprf
To see the claim regarding composition, let $I\xrightarrow{\sigma}\Alg{A}\otimes\Alg{E}$ be a dilation of $\alpha$. Then 
\begingroup
\allowdisplaybreaks
\begin{align*}
\vcenter{\hbox{%
\begin{tikzpicture}[font=\small]
\node[state] (omega) at (0,0) {\;\;$\sigma$\;\;};
\coordinate (X) at (0.45,2.65) {};
\coordinate (Y) at (-0.45,2.65) {};
\node[arrow box] (psi) at (-0.45,0.85) {\;$\map{E}_{1}$\;};
\node[arrow box] (e) at (-0.45,1.85) {\;$\map{F}_{1}$\;};
\draw (omega) ++(0.45, 0) to (X);
\draw (omega) ++(-0.45, 0) to (psi);
\draw (psi) to (e);
\draw (e) to (Y);
\path[scriptstyle]
node at (0.65,2.45) {$\Alg{E}$}
node at (-0.65,1.35) {$\Alg{B}$}
node at (-0.65,2.45) {$\Alg{C}$}
node at (-0.65,0.3) {$\Alg{A}$};
\end{tikzpicture}}}
\;\;&=\;\;
\vcenter{\hbox{%
\begin{tikzpicture}[font=\small]
\node[state] (omega) at (0,0) {\;\;$\sigma$\;\;};
\coordinate (X) at (0.45,2.65) {};
\coordinate (Y) at (-0.45,2.65) {};
\node[arrow box] (psi) at (-0.45,0.85) {\;$\map{E}_{2}$\;};
\node[arrow box] (e) at (-0.45,1.85) {\;$\map{F}_{1}$\;};
\draw (omega) ++(0.45, 0) to (X);
\draw (omega) ++(-0.45, 0) to (psi);
\draw (psi) to (e);
\draw (e) to (Y);
\path[scriptstyle]
node at (0.65,2.45) {$\Alg{E}$}
node at (-0.65,1.35) {$\Alg{B}$}
node at (-0.65,2.45) {$\Alg{C}$}
node at (-0.65,0.3) {$\Alg{A}$};
\end{tikzpicture}}}
&&\mbox{because $\map{E}_{1}\aeequals{\alpha}\map{E}_{2}$ and $\sigma$ is a dilation of $\alpha$}\\
\;\;&=\;\;
\vcenter{\hbox{%
\begin{tikzpicture}[font=\small]
\node[state] (omega) at (0,0) {\;\;$\sigma$\;\;};
\coordinate (X) at (0.45,2.65) {};
\coordinate (Y) at (-0.45,2.65) {};
\node[arrow box] (psi) at (-0.45,0.85) {\;$\map{E}_{2}$\;};
\node[arrow box] (e) at (-0.45,1.85) {\;$\map{F}_{2}$\;};
\draw (omega) ++(0.45, 0) to (X);
\draw (omega) ++(-0.45, 0) to (psi);
\draw (psi) to (e);
\draw (e) to (Y);
\path[scriptstyle]
node at (0.65,2.45) {$\Alg{E}$}
node at (-0.65,1.35) {$\Alg{B}$}
node at (-0.65,2.45) {$\Alg{C}$}
node at (-0.65,0.3) {$\Alg{A}$};
\end{tikzpicture}}}
&&\mbox{because $\map{F}_{1}\aeequals{\beta}\map{F}_{2}$ and $(\map{E}_{2}\otimes\id_{\Alg{E}})\circ\sigma$ is a dilation of $\beta$,}\\
\end{align*}
\endgroup
which says 
\[
\begin{split}
\big((\map{F}_{1}\circ\map{E}_{1})\otimes\id_{\Alg{E}}\big)\circ\sigma
&=(\map{F}_{1}\otimes\id_{\Alg{E}})\circ\big((\map{E}_{1}\otimes\id_{\Alg{E}})\circ\sigma\big)
=(\map{F}_{1}\otimes\id_{\Alg{E}})\circ\big((\map{E}_{2}\otimes\id_{\Alg{E}})\circ\sigma\big)\\
&=(\map{F}_{2}\otimes\id_{\Alg{E}})\circ\big((\map{E}_{2}\otimes\id_{\Alg{E}})\circ\sigma\big)
=\big((\map{F}_{2}\circ\map{E}_{2})\otimes\id_{\Alg{E}}\big)\circ\sigma,
\end{split}
\]
where parentheses have been included for emphasis. 
To see the claim regarding tensor products, let $I\xrightarrow{\tau}\Alg{A}\otimes\Alg{A}'\otimes\Alg{F}$ be a dilation of $\alpha\otimes\alpha'$. 
A proof can be obtained via string diagrams as follows: 
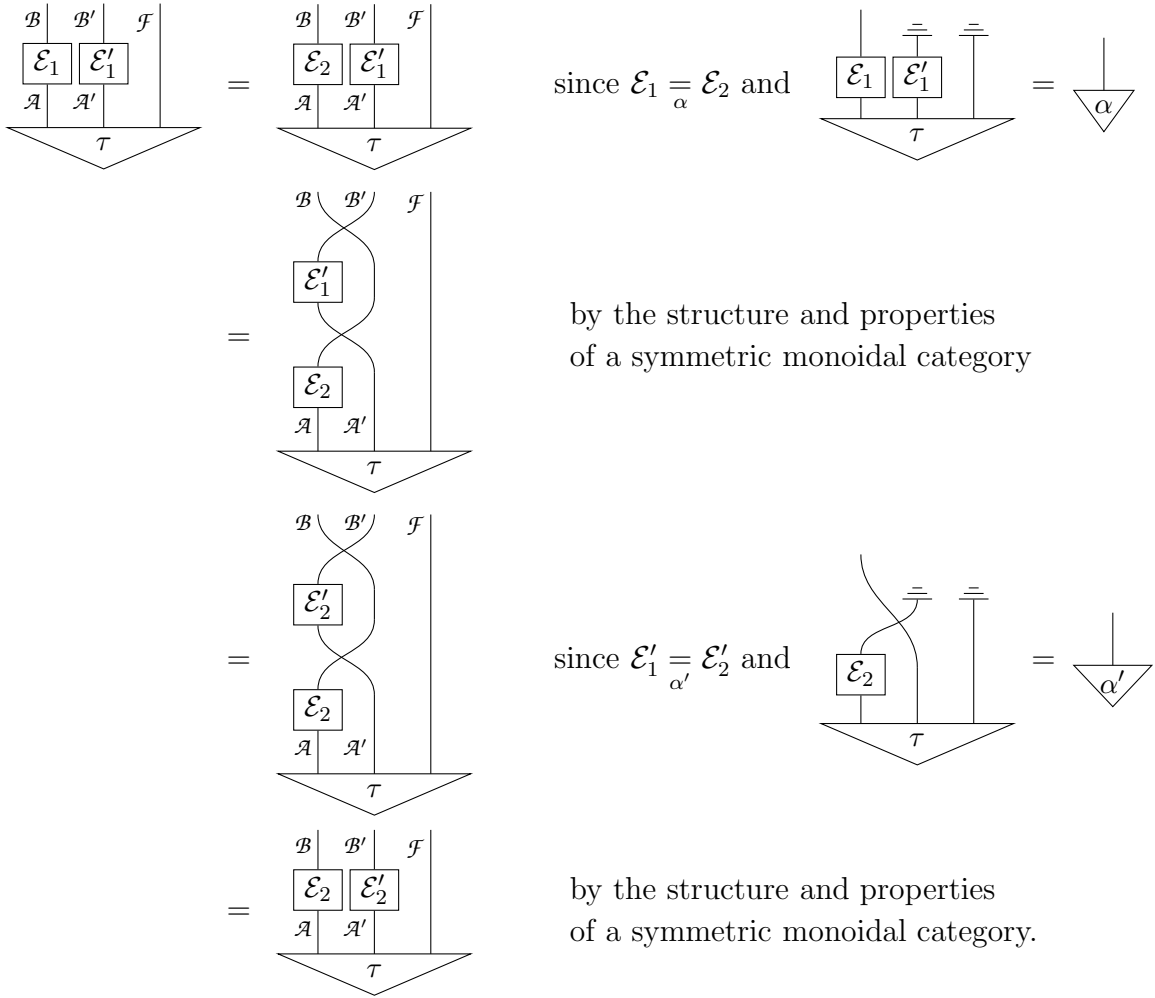
\begingroup
\allowdisplaybreaks
\begin{align*}
\vcenter{\hbox{%
\begin{tikzpicture}[font=\small]
\node[state] (omega) at (0,0) {\;\;\;$\tau$\;\;\;};
\coordinate (X) at (-0.75,1.65) {};
\coordinate (Y) at (0.0,1.65) {};
\coordinate (Z) at (0.75,1.65) {};
\node[arrow box] (h) at (-0.75,0.85) {$\map{E}_{1}$};
\node[arrow box] (g) at (0.0,0.85) {$\map{E}_{1}'$};
\draw (omega) ++(-0.75, 0) to (h);
\draw (omega) ++(0.0, 0) to (g);
\draw (h) to (X);
\draw (g) to (Y);
\draw (omega) ++(0.75, 0) to (Z);
\path[scriptstyle]
node at (-0.95,1.45) {$\Alg{B}$}
node at (-0.95,0.35) {$\Alg{A}$}
node at (-0.25,1.45) {$\Alg{B}'$}
node at (-0.25,0.35) {$\Alg{A}'$}
node at (0.55,1.40) {$\Alg{F}$};
\end{tikzpicture}}}
\;\;&=\;\;
\vcenter{\hbox{%
\begin{tikzpicture}[font=\small]
\node[state] (omega) at (0,0) {\;\;\;$\tau$\;\;\;};
\coordinate (X) at (-0.75,1.65) {};
\coordinate (Y) at (0.0,1.65) {};
\coordinate (Z) at (0.75,1.65) {};
\node[arrow box] (h) at (-0.75,0.85) {$\map{E}_{2}$};
\node[arrow box] (g) at (0.0,0.85) {$\map{E}_{1}'$};
\draw (omega) ++(-0.75, 0) to (h);
\draw (omega) ++(0.0, 0) to (g);
\draw (h) to (X);
\draw (g) to (Y);
\draw (omega) ++(0.75, 0) to (Z);
\path[scriptstyle]
node at (-0.95,1.45) {$\Alg{B}$}
node at (-0.95,0.35) {$\Alg{A}$}
node at (-0.25,1.45) {$\Alg{B}'$}
node at (-0.25,0.35) {$\Alg{A}'$}
node at (0.55,1.45) {$\Alg{F}$};
\end{tikzpicture}}}
&&\mbox{since $\map{E}_{1}\aeequals{\alpha}\map{E}_{2}$ and \; $\vcenter{\hbox{%
\begin{tikzpicture}[font=\small]
\node[state] (omega) at (0,0) {\;\;\;$\tau$\;\;\;};
\coordinate (X) at (-0.75,1.45) {};
\node[discarder] (Y) at (0.0,1.10) {};
\node[discarder] (Z) at (0.75,1.10) {};
\node[arrow box] (h) at (-0.75,0.55) {$\map{E}_{1}$};
\node[arrow box] (g) at (0.0,0.55) {$\map{E}_{1}'$};
\draw (omega) ++(-0.75, 0) to (h);
\draw (omega) ++(0.0, 0) to (g);
\draw (h) to (X);
\draw (g) to (Y);
\draw (omega) ++(0.75, 0) to (Z);
\end{tikzpicture}}}
\;=\;
\vcenter{\hbox{%
\begin{tikzpicture}[font=\small]
\node[state] (omega) at (0,-0.8) {$\alpha$};
\coordinate (q) at (0,-0.1);
\draw (omega) to (q);
\end{tikzpicture}}}
$ 
}
\\
\;\;&=\;\;\vcenter{\hbox{%
\begin{tikzpicture}[font=\small]
\node[state] (omega) at (0,0) {\;\;\;$\tau$\;\;\;};
\coordinate (X) at (-0.75,3.45) {};
\coordinate (Y) at (0.0,3.45) {};
\coordinate (Z) at (0.75,3.45) {};
\node[arrow box] (h) at (-0.75,0.85) {$\map{E}_{2}$};
\coordinate (2l) at (0,1.05);
\coordinate (2ma) at (0,2.05);
\coordinate (2mb) at (0,2.45);
\node[arrow box] (g) at (-0.75,2.25) {$\map{E}_{1}'$};
\draw (omega) ++(-0.75, 0) to (h);
\draw (omega) ++(0.0, 0) to (2l);
\draw (2l) to[out=90,in=-90] (g);
\draw (h) to[out=90,in=-90] (2ma);
\draw (2ma) to (2mb); 
\draw (2mb) to[out=90,in=-90] (X);
\draw (g) to[out=90,in=-90] (Y);
\draw (omega) ++(0.75, 0) to (Z);
\path[scriptstyle]
node at (-0.95,3.35) {$\Alg{B}$}
node at (-0.95,0.35) {$\Alg{A}$}
node at (-0.25,3.35) {$\Alg{B}'$}
node at (-0.25,0.35) {$\Alg{A}'$}
node at (0.55,3.30) {$\Alg{F}$};
\end{tikzpicture}}}
&&\mbox{\begin{tabular}{l}by the structure and properties\\of a symmetric monoidal category\end{tabular}}
\\
\;\;&=\;\;\vcenter{\hbox{%
\begin{tikzpicture}[font=\small]
\node[state] (omega) at (0,0) {\;\;\;$\tau$\;\;\;};
\coordinate (X) at (-0.75,3.45) {};
\coordinate (Y) at (0.0,3.45) {};
\coordinate (Z) at (0.75,3.45) {};
\node[arrow box] (h) at (-0.75,0.85) {$\map{E}_{2}$};
\coordinate (2l) at (0,1.05);
\coordinate (2ma) at (0,2.05);
\coordinate (2mb) at (0,2.45);
\node[arrow box] (g) at (-0.75,2.25) {$\map{E}_{2}'$};
\draw (omega) ++(-0.75, 0) to (h);
\draw (omega) ++(0.0, 0) to (2l);
\draw (2l) to[out=90,in=-90] (g);
\draw (h) to[out=90,in=-90] (2ma);
\draw (2ma) to (2mb); 
\draw (2mb) to[out=90,in=-90] (X);
\draw (g) to[out=90,in=-90] (Y);
\draw (omega) ++(0.75, 0) to (Z);
\path[scriptstyle]
node at (-0.95,3.35) {$\Alg{B}$}
node at (-0.95,0.35) {$\Alg{A}$}
node at (-0.25,3.35) {$\Alg{B}'$}
node at (-0.25,0.35) {$\Alg{A}'$}
node at (0.55,3.30) {$\Alg{F}$};
\end{tikzpicture}}}
&&\mbox{since $\map{E}_{1}'\aeequals{\alpha'}\map{E}_{2}'$ and \; 
$\vcenter{\hbox{%
\begin{tikzpicture}[font=\small]
\node[state] (omega) at (0,0) {\;\;\;$\tau$\;\;\;};
\coordinate (X) at (-0.75,2.05);
\node[discarder] (Z) at (0.75,1.65) {};
\node[arrow box] (h) at (-0.75,0.65) {$\map{E}_{2}$};
\coordinate (2l) at (0,0.75);
\node[discarder] (2ma) at (0,1.65) {};
\coordinate (g) at (-0.75,2.25);
\draw (omega) ++(-0.75, 0) to (h);
\draw (omega) ++(0.0, 0) to (2l);
\draw (2l) to[out=90,in=-90] (g);
\draw (h) to[out=90,in=-90] (2ma);
\draw (omega) ++(0.75, 0) to (Z);
\end{tikzpicture}}}
\;=\;
\vcenter{\hbox{%
\begin{tikzpicture}[font=\small]
\node[state] (omega) at (0,-0.8) {$\alpha'$};
\coordinate (q) at (0,-0.1);
\draw (omega) to (q);
\end{tikzpicture}}}
$
}
\\
\;\;&=\;\;
\vcenter{\hbox{%
\begin{tikzpicture}[font=\small]
\node[state] (omega) at (0,0) {\;\;\;$\tau$\;\;\;};
\coordinate (X) at (-0.75,1.65) {};
\coordinate (Y) at (0.0,1.65) {};
\coordinate (Z) at (0.75,1.65) {};
\node[arrow box] (h) at (-0.75,0.85) {$\map{E}_{2}$};
\node[arrow box] (g) at (0.0,0.85) {$\map{E}_{2}'$};
\draw (omega) ++(-0.75, 0) to (h);
\draw (omega) ++(0.0, 0) to (g);
\draw (h) to (X);
\draw (g) to (Y);
\draw (omega) ++(0.75, 0) to (Z);
\path[scriptstyle]
node at (-0.95,1.45) {$\Alg{B}$}
node at (-0.95,0.35) {$\Alg{A}$}
node at (-0.25,1.45) {$\Alg{B}'$}
node at (-0.25,0.35) {$\Alg{A}'$}
node at (0.55,1.40) {$\Alg{F}$};
\end{tikzpicture}}}
&&\mbox{\begin{tabular}{l}by the structure and properties\\of a symmetric monoidal category.\end{tabular}}
\end{align*}
\endgroup
This proves that $\map{E}_{1}\otimes\map{E}_{1}'\aeequals{\alpha\otimes\alpha'}\map{E}_{2}\otimes\map{E}_{2}'$, thus concluding the proof. 
\eprf

\bn
\label{prop:statesandaemorphisms}
Let $\mathbf{C}$ be a semicartesian category. 
The collection of objects given by pairs $(\Alg{A},\alpha)$ with $\Alg{A}$ in $\mathbf{C}$ together with the collection of dilationally equal equivalence classes of state-preserving morphisms $(\Alg{A},\alpha)\xrightarrow{\map{E}}(\Alg{B},\beta)$ defines a symmetric monoidal (in fact, semicartesian) category, denoted by $\states(\mathbf{C})$. In addition, $\states^{+}(\mathbf{C})$ embeds fully and faithfully into $\states(\mathbf{C})$. 
\en

\bprf
The well-definedness of composition and the tensor product follows from Lemma~\ref{lem:dilationcompandtensor}. The full and faithful embedding follows from the definitions of dilational equality and faithful states.
\eprf

\section{Categorical retrodiction and probability kinematics}
\label{sec:retrodiction}

\bd
Let $\mathbf{C}$ be a semicartesian category. 
Using the same notation as in Definition~\ref{defn:fstatesSCC} and Proposition~\ref{prop:statesandaemorphisms}, a \define{retrodiction family on faithful states} for $\mathbf{C}$ is an assignment $\retro:\states^{+}(\mathbf{C})\to\states^{+}(\mathbf{C})^{\op}$ that acts as the identity on objects. 
More generally, a \define{retrodiction family} for $\mathbf{C}$ is an assignment $\retro:\states(\mathbf{C})\to\states(\mathbf{C})^{\op}$ that acts as the identity on objects. 
Given a morphism $(\Alg{A},\alpha)\xrightarrow{\map{E}}(\Alg{B},\beta)$ in $\states^{+}(\mathbf{C})$ (or $\states(\mathbf{C})$), the morphism $(\Alg{B},\beta)\xrightarrow{\retro_{\alpha,\map{E}}:=\retro(\map{E})}(\Alg{A},\alpha)$ is called a \define{recovery morphism} for $(\alpha,\map{E})$. The map $\retro_{\alpha,\map{E}}\in\states(\mathbf{C})^{\op}$ (or a \emph{representative} for $\retro_{\alpha,\map{E}}\in\states^{+}(\mathbf{C})^{\op})$ is called a \define{recovery map}. 
A retrodiction family is thus said to satisfy the \define{recovery property}, i.e., is \define{recovering}. 
\ed

\br
The terminology of \emph{recovery map} is used to be consistent with the quantum information literature~\cite{LiWi18,SuToHa16,CHPSSW19,JRSWW16,JRSWW18}.  
This is because $\retro_{\alpha,\map{E}}\circ\beta=\alpha$, a condition that is typically referred to as the \emph{recovery property}, follows from two conditions that the objects are fixed together with the fact that $\retro$ is contravariant. 
Note that this does \emph{not} imply $\retro_{\alpha,\map{E}}\circ\map{E}=\id_{\Alg{A}}$ (or even up to dilational equality), which is a condition closely related to perfect error-correcting codes and disintegrations~\cite{PaRu19,PaBayes}. Here, $\retro_{\alpha,\map{E}}$ acts more as a state-specific approximate error-correcting code~\cite{BaKn02,NgMa10}. 
The word \emph{family} in \emph{retrodiction family} is used to avoid calling it a \emph{functor} yet, since there are many useful retrodiction families in the literature that do not satisfy functoriality~\cite{PaBu22}. 
Meanwhile, the notion of a \emph{retrodiction family} on faithful states was introduced for the categories $\FinStoch$ and $\NCFinStoch$ in Ref.~\cite{PaBu22}. 
In the present paper, we have extended this definition to arbitrary semicartesian categories and also states that are not necessarily faithful. 
\er

\bd
Using the same notation as in Definition~\ref{defn:fstatesSCC} and Proposition~\ref{prop:statesandaemorphisms}, either let $\mathbf{C}$ be a semicartesian category for which $\states^{+}(\mathbf{C})$ is a symmetric monoidal category or let $\mathbf{C}$ be any semicartesian category, so that $\states(\mathbf{C})$ is automatically a symmetric monoidal category. 
A \define{retrodiction functor on faithful states} for $\mathbf{C}$ is an assignment $\retro:\states^{+}(\mathbf{C})\to\states^{+}(\mathbf{C})^{\op}$ that is a retrodiction family for faithful states and an inverting monoidal involutive functor.
A \define{retrodiction functor} for $\mathbf{C}$ is an assignment $\retro:\states(\mathbf{C})\to\states(\mathbf{C})^{\op}$ that is a retrodiction family and an inverting monoidal involutive functor.
Explicitly, the conditions state that $\retro$ is 
\begin{enumerate}
\item
\define{recovering} in the sense that $\retro$ is a retrodiction family, i.e., $\retro$ fixes objects.
\item
\define{normalizing} in the sense that retrodicting the process that does nothing is also the process that does nothing, i.e., $\retro_{\alpha,\id_{\Alg{A}}}=\id_{\Alg{A}}$. 
\item
\define{compositional} in the sense that if $(\Alg{A},\alpha)\xrightarrow{\map{E}}(\Alg{B},\beta)\xrightarrow{\map{F}}(\Alg{C},\gamma)$ is a composable pair of morphisms, then the retrodiction $\retro_{\alpha,\map{F}\circ\map{E}}$ associated with the composite process $\map{F}\circ\map{E}$ is the composite of the retrodictions $\retro_{\map{E}(\alpha),\map{F}}:\Alg{C}\to\Alg{B}$ and $\retro_{\alpha,\map{E}}:\Alg{B}\to\Alg{A}$ associated with the constituent components, i.e., $\retro_{\alpha,\map{F}\circ\map{E}}=\retro_{\alpha,\map{E}}\circ\retro_{\map{E}(\alpha),\map{F}}$. 
\item
\define{tensorial} in the sense given two pairs $(\Alg{A},\alpha)\xrightarrow{\map{E}}(\Alg{B},\beta)$ and $(\Alg{A}',\alpha')\xrightarrow{\map{E}'}(\Alg{B}',\beta')$, then the retrodiction $\retro_{\alpha\otimes\alpha',\map{E}\otimes\map{E}'}:\Alg{B}\otimes\Alg{B}'\to\Alg{A}\otimes\Alg{A}'$ associated to the tensor product of the systems and processes is equal to the tensor product $\retro_{\alpha,\map{E}}\otimes\retro_{\alpha',\map{E}'}$ of the constituent retrodictions, i.e., $\retro_{\alpha\otimes\alpha',\map{E}\otimes\map{E}'}=\retro_{\alpha,\map{E}}\otimes\retro_{\alpha',\map{E}'}$.
\item
\define{inverting} in the sense that the retrodiction of an isomorphism $(\Alg{A},\alpha)\xrightarrow{\map{E}}(\Alg{B},\beta)$ is the inverse $\map{E}^{-1}$ of the original process, i.e., $\retro_{\alpha,\map{E}}=\map{E}^{-1}$.
\item
\define{involutive} in the sense that retrodicting a retrodiction gives back the original process, i.e., $\retro_{\map{E}(\alpha),\retro_{\alpha,\map{E}}}=\map{E}$, or more cleanly $\map{R}\circ\map{R}=\id_{\states^{+}(\mathbf{C})}$ (or $\map{R}\circ\map{R}=\id_{\states(\mathbf{C})}$).
\end{enumerate}
\ed

The normalizing and compositional axioms say that $\retro$ is a functor. Adding the tensorial axiom says that $\retro$ is a monoidal functor. Adding the recovery and involutive axioms says that $\retro$ is a monoidal dagger. Finally, the addition of the inverting condition says that $\retro$ is an inverting monoidal dagger. Before giving examples, we first explain why we use a category of states on $\mathbf{C}$ and state-preserving processes rather than $\mathbf{C}$ itself. We then explain how having a retrodiction functor enables a generalization of Jeffrey's probability kinematics and Bayesian updating~\cite{Ja19,Je90}. 

\br
There does not exist a recovering monoidal inverting involutive functor on $\FinStoch$ or $\NCFinStoch$. Indeed, if such a functor $\retro$ existed (in fact, if it was just assumed that $\retro$ was recovering and involutive), then given any two different states $I\xrightarrow{\alpha_{1}}\Alg{A}$ and $I\xrightarrow{\alpha_{2}}\Alg{A}$, one has $\retro(\alpha_{1})=\;!_{\Alg{A}}=\retro(\alpha_{2})$. Applying $\retro$ again and invoking involutivity would lead to the contradiction $\alpha_{1}=\alpha_{2}$. Such observations have led many to believe that such an involutive functor must therefore only exist for a suitable proper subcategory of $\FinStoch$ or $\NCFinStoch$~\cite{CGS17,CAZ21,ChLi22}, which is often taken to be the subcategory of bistochastic maps or unital quantum channels, respectively. 
However, these subcategories are \emph{isomorphic} to the full subcategories of $\states^{+}(\FinStoch)$ and $\states^{+}(\NCFinStoch)$ whose objects are of the form $(\Alg{A},u_{\Alg{A}})$, where $u_{\Alg{A}}:=\frac{1_{\Alg{A}}}{\Tr[1_{\Alg{A}}]}$ denotes the \define{uniform} (or \define{maximally mixed}) state on $\Alg{A}$ (and $\Alg{A}$ is commutative in the case of $\FinStoch$). 
This observation led Ref.~\cite{PaBu22} to propose to use a category of states and state-preserving processes to model a category of channels that admits a  recovering monoidal inverting involutive functor on it, thus dispelling the belief that such a form of time-reversal symmetry is not possible for all quantum channels~\cite{CGS17,CAZ21,ChLi22}. See Ref.~\cite{PaBu22} for more details. 
\er

\bd
Let $\mathbf{C}$ be a semicartesian category. A state $I\xrightarrow{\alpha}\Alg{A}$ in $\mathbf{C}$ is said to be \define{absolutely continuous} with respect to a state $I\xrightarrow{\tilde{\alpha}}\Alg{A}$, written $\alpha\ll\tilde{\alpha}$ or $\tilde{\alpha}\gg\alpha$, iff for every pair of morphisms $\map{E},\map{F}:\Alg{A}\to\Alg{B}$ such that $\map{E}\aeequals{\tilde{\alpha}}\map{F}$, then $\map{E}\aeequals{\alpha}\map{F}$. 
\ed

\bx
By Proposition~\ref{prop:dilationalequalityQuantum}, in the case of $\NCFinStoch$, a state $\C\xrightarrow{\alpha}\Alg{A}$ is absolutely continuous with respect to a state $\C\xrightarrow{\tilde{\alpha}}\Alg{A}$ whenever $P_{\alpha}\le P_{\tilde{\alpha}}$ for the support projections (the $\le$ here refers to the operator partial ordering in terms of positive elements~\cite{Bh07,Fa01}). This agrees with the standard definition of absolute continuity in the noncommutative setting~\cite{Dy52,Hi84}.
\ex

\begin{term}
Let $\mathbf{C}$ be a semicartesian category. 
Given a morphism $(\Alg{A},\alpha)\xrightarrow{\map{E}}(\Alg{B},\beta)$ in $\states(\mathbf{C})$, if $I\xrightarrow{\epsilon}\Alg{B}$ is a state satisfying $\epsilon\ll\beta$, called \define{soft evidence}, then 
$\retro_{\alpha,\map{E}}\circ\epsilon$
is the \define{updated} state on $\Alg{A}$ based on the evidence $\epsilon$. This procedure of pushing evidence $\epsilon$ backwards along $\retro_{\alpha,\map{E}}$ is called \define{Jeffrey's probability kinematics} and generalizes Bayes' rule. One can visualize Jeffrey's probability kinematics schematically via 
\[
\xy0;/r.25pc/:
(0,6)*+{I}="1";
(-6,-6)*+{\Alg{A}}="X";
(6,-6)*+{\Alg{B}}="Y";
{\ar"1";"X"_{\alpha}};
{\ar"1";"Y"^{\beta}};
{\ar"X";"Y"_{\map{E}}};
\endxy
\quad\text{ and }\quad
\xy0;/r.25pc/:
(0,6)*+{I}="1";
(6,-6)*+{\Alg{B}}="Y";
{\ar"1";"Y"^{\epsilon}};
\endxy
\ll
\xy0;/r.25pc/:
(0,6)*+{I}="1";
(6,-6)*+{\Alg{B}}="Y";
{\ar"1";"Y"^{\beta}};
\endxy
\quad\text{ give }\quad
\xy0;/r.25pc/:
(0,6)*+{I}="1";
(-6,-6)*+{\Alg{A}}="X";
(6,-6)*+{\Alg{B}}="Y";
{\ar"1";"X"_{\retro_{\alpha,\map{E}}\circ\epsilon}};
{\ar"1";"Y"^{\epsilon}};
{\ar"Y";"X"^{\retro_{\alpha,\map{E}}}};
\endxy
\]
\end{term}

\br
The condition $\epsilon\ll\beta$ guarantees that $\retro_{\alpha,\map{E}}\circ\epsilon$ is well-defined, i.e., independent of the choice of representative $\retro_{\alpha,\map{E}}$. This is because if $\retro_{\alpha,\map{E}}'$ is another representative, then $\retro_{\alpha,\map{E}}\aeequals{\beta}\retro_{\alpha,\map{E}}'$ by definition. The assumption $\epsilon\ll\beta$ then guarantees $\retro_{\alpha,\map{E}}\aeequals{\epsilon}\retro_{\alpha,\map{E}}'$, which implies $\retro_{\alpha,\map{E}}\circ\epsilon=\retro_{\alpha,\map{E}}'\circ\epsilon$. If one restricts to the category $\states^{+}(\mathbf{C})$ of faithful states and state-preserving morphisms, then one can ignore the absolutely continuous conditions on $\epsilon$ since they automatically hold. 
\er

\bx
\label{ex:FinStochretro}
We have already explained how Bayesian inversion defines a retrodiction functor for $\FinStoch$ in Theorem~\ref{thm:BisRetro}. The observation that Bayesian inversion defines a monoidal dagger (the assumptions of which include most of the properties of a retrodiction functor) seems to have first been made in the more general setting of Borel spaces in~\cite[Theorem 2.10]{DSDG18} and in a more general categorical setting of Markov categories with conditionals in~\cite[Remark~13.10]{Fr20}. Together with the results of Ref.~\cite{PaBayes} regarding the inverting property of categorical Bayesian inverses, this shows that Markov categories with conditionals exhibit retrodiction functors. 
\ex

\bq
Does every Markov category with conditionals admit a \emph{unique} retrodiction functor? Note that an affirmative answer to this question does not guarantee the uniqueness of a retrodiction functor for $\NCFinStoch$ because the latter is not a Markov category due to the no-broadcasting theorem (and even if one were to view $\NCFinStoch$ as embedded in a quantum Markov category, $\NCFinStoch$ would not admit conditionals)~\cite{PaBayes,Um54}. 
\eq

\bx
\label{defn:Petzmaps}
Let $(\Alg{A},\alpha)\xrightarrow{\map{E}}(\Alg{B},\beta)$ be a morphism in $\states^{+}(\NCFinStoch)$. The \define{Petz recovery map} associated with the pair $(\alpha,\map{E})$ is the morphism
\[
(\Alg{A},\alpha)\xleftarrow{\retro^{\Petz}_{\alpha,\mathcal{E}}}(\Alg{B},\beta)
\]
in $\states^{+}(\NCFinStoch)$ defined by the formula
\[
\retro^{\Petz}_{\alpha,\map{E}}
:=\Ad_{\alpha^{1/2}}\circ \map{E}^*\circ\Ad_{\beta^{-1/2}}
,
\]
where $\map{E}^*$ is the Hilbert--Schmidt adjoint of $\map{E}$ and $\Ad_{V}$ is the map sending $A$ to $VAV^{\dag}$ (see Refs.~\cite{Wilde15,PaBu22} for details). Then $\retro^{\Petz}:\states^{+}(\NCFinStoch)\to\states^{+}(\NCFinStoch)^{\op}$ defines a retrodiction functor on $\NCFinStoch$~\cite{LiWi18,Wilde15,PaBu22}. 
\ex

More generally, we can define a retrodiction functor for all states, not necessarily faithful. 

\bt
\label{thm:Petzretrodictionae}
The assignment 
\[
\begin{split}
\states(\NCFinStoch)&\xrightarrow{\retro^{\Petz}}\states(\NCFinStoch)^{\op}\\
(\Alg{A},\alpha)\xrightarrow{[\map{E}]}(\Alg{B},\beta)&\mapsto(\Alg{A},\alpha)\xleftarrow{[\retro^{\Petz}_{\alpha,\map{E}}]}(\Alg{B},\beta), 
\end{split}
\]
where 
\[
\Alg{B}\ni B\mapsto\retro^{\Petz}_{\alpha,\map{E}}(B)
:=
(\Ad_{\alpha^{1/2}}\circ\map{E}^*\circ\Ad_{\hat{\beta}^{1/2}})(B)+ \Tr(P_{\beta}^{\perp}B)\frac{P_{\alpha}}{\Tr[P_{\alpha}]} 
\]
defines a retrodiction functor on $\NCFinStoch$. Here, $\hat{\beta}$ denotes the Moore--Penrose pseudo-inverse of $\beta$~\cite{Mo1920,Pe55} and the $[\;\cdot\;]$ bracket denotes a dilationally equal equivalence class with respect to the state on the domain. 
\et

Although Ref.~\cite{Wilde15} gives a proof of this, we provide additional details to justify certain steps (specifically the compositionality property). We will utilize one lemma regarding support projections. 

\blem
\label{lem:nullspaceE}
Let $(\Alg{A},\alpha)\xrightarrow{\map{E}}(\Alg{B},\beta)$ be a state-preserving CPTP map. Then $\Ad_{P_{\alpha}}\circ\map{E}^*\circ\Ad_{P_{\beta}}=\Ad_{P_{\alpha}}\circ\map{E}^*$ as linear maps from $\Alg{B}$ to $\Alg{A}$. 
\elem

\bprf[Proof of Lemma~\ref{lem:nullspaceE}]
Let 
\[
\Alg{N}_{\,\alpha}:=\big\{A\in\Alg{A}\;:\;\Tr[\alpha A^{\dag}A]=0\big\}=\big\{AP_{\alpha}^{\perp}\;:\;A\in\Alg{A}\big\}=:\Alg{A}P_{\alpha}^{\perp}
\]
denote the \define{nullspace} of $\alpha$, and similarly for $\beta$~\cite{GN43,Se47,PaGNS}. Then $\map{E}^*(\Alg{N}_{\,\beta})\subseteq\Alg{N}_{\,\alpha}$ by~\cite[Lemma 2.31]{GPRR21}. This means that $\map{E}^*(BP_{\beta}^{\perp})\in\Alg{A}P_{\alpha}^{\perp}$ for all $B\in\Alg{B}$. It also means that $\map{E}^*(P_{\beta}^{\perp}B)=\big(\map{E}^*(B^{\dag}P_{\beta}^{\perp})\big)^{\dag}\in P_{\alpha}^{\perp}\Alg{A}$ since $\map{E}^*$ is dagger-preserving, i.e., $\map{E}^*(B^{\dag})^{\dag}=\map{E}^*(B)$ for all $B\in\Alg{B}$. All these facts imply 
\[
\Ad_{P_{\alpha}}\big(\map{E}^*(B)\big)
=\Ad_{P_{\alpha}}\big(\map{E}^*(P_{\beta}BP_{\beta}+P_{\beta}BP_{\beta}^{\perp}+P_{\beta}^{\perp}BP_{\beta}+P_{\beta}^{\perp}BP_{\beta}^{\perp})\big)
=\Ad_{P_{\alpha}}\Big(\map{E}^*\big(\Ad_{P_{\beta}}(B)\big)\Big)
\]
for all $B\in\Alg{B}$. This proves the claim. 
\eprf

\bprf[Proof of Theorem~\ref{thm:Petzretrodictionae}]
To see that $\retro^{\Petz}$ is well-defined, suppose $\map{E}_{1}\aeequals{\alpha}\map{E}_{2}$ for two state-preserving morphisms $\map{E}_{1},\map{E}_{2}:(\Alg{A},\alpha)\to(\Alg{B},\beta)$. Then
\begin{align*}
\retro^{\Petz}_{\alpha,\map{E}_{1}}\circ\Ad_{P_{\beta}}
&=\Ad_{\alpha^{1/2}}\circ\Ad_{P_{\alpha}}\circ\map{E}_{1}^*\circ\Ad_{\hat{\beta}^{1/2}} &&\mbox{since $\hat{\beta}^{1/2}P_{\beta}=\hat{\beta}^{1/2}$, $P_{\beta}^{\perp}P_{\beta}=0$, and $\alpha^{1/2}P_{\alpha}=\alpha^{1/2}$}\\
&=\Ad_{\alpha^{1/2}}\circ\big(\map{E}_{1}\circ\Ad_{P_{\alpha}}\big)^{*}\circ\Ad_{\hat{\beta}^{1/2}} &&\mbox{since $P_{\alpha}^{\dag}=P_{\alpha}$ and $(\map{F}\circ\map{E})^*=\map{E}^*\circ\map{F}^*$} \\
&=Ad_{\alpha^{1/2}}\circ\big(\map{E}_{2}\circ\Ad_{P_{\alpha}}\big)^{*}\circ\Ad_{\hat{\beta}^{1/2}} &&\mbox{by Proposition~\ref{prop:dilationalequalityQuantum} and because $\map{E}_{1}\aeequals{\alpha}\map{E}_{2}$}\\
&=\retro^{\Petz}_{\alpha,\map{E}_{2}}\circ\Ad_{P_{\beta}}&&\mbox{by reversing the calculation.}
\end{align*}
By Proposition~\ref{prop:dilationalequalityQuantum}, this shows that $\retro^{\Petz}_{\alpha,\map{E}_{1}}\aeequals{\beta}\retro^{\Petz}_{\alpha,\map{E}_{2}}$, which proves that $\retro^{\Petz}$ is well-defined.  
The fact that $\retro^{\Petz}_{\alpha,\map{E}}$ defines a state-preserving CPTP map is proved in Ref.~\cite{Wilde15}. 
The fact that $\retro^{\Petz}$ is normalizing and tensorial is proved in Ref.~\cite{Wilde15}. 
The fact that $\retro^{\Petz}$ is inverting is proved in Refs.~\cite{BuSc21,PaBu22}.

To see that $\retro^{\Petz}$ is compositional, let $(\Alg{A},\alpha)\xrightarrow{\map{E}}(\Alg{B},\beta)\xrightarrow{\map{F}}(\Alg{C},\gamma)$ be a pair of composable state-preserving CPTP maps. Then
\begin{align*}
\retro^{\Petz}_{\alpha,\map{F}\circ\map{E}}\circ\Ad_{P_{\gamma}}
&=\Ad_{\alpha^{1/2}}\circ(\map{F}\circ\map{E})^*\circ\Ad_{\hat{\gamma}^{1/2}}&& \mbox{by definition}\\
&=\Ad_{\alpha^{1/2}}\circ\map{E}^*\circ\map{F}^*\circ\Ad_{\hat{\gamma}^{1/2}} && \mbox{by a property of the Hilbert--Schmidt adjoint}\\
&=\Ad_{\alpha^{1/2}}\circ\map{E}^*\circ\Ad_{P_{\beta}}\circ\map{F}^*\circ\Ad_{\hat{\gamma}^{1/2}} && \mbox{by Lemma~\ref{lem:nullspaceE} and $\Ad_{\alpha^{1/2}}=\Ad_{\alpha^{1/2}}\circ\Ad_{P_{\alpha}}$}\\
&=\retro^{\Petz}_{\alpha,\map{E}}\circ\retro^{\Petz}_{\beta,\map{F}}\circ\Ad_{P_{\gamma}} &&\mbox{by direct calculation and since $\hat{\beta}\beta=P_{\beta}$.}
\end{align*}
Proposition~\ref{prop:dilationalequalityQuantum} then shows $\retro^{\Petz}_{\alpha,\map{E}}\circ\retro^{\Petz}_{\beta,\map{F}}\aeequals{\gamma}\retro^{\Petz}_{\alpha,\map{F}\circ\map{E}}$.

The final thing to show is that $\retro^{\Petz}$ is involutive, i.e., $\retro^{\Petz}_{\beta,\retro^{\Petz}_{\alpha,\map{E}}}\aeequals{\alpha}\map{E}$ for all $(\Alg{A},\alpha)\xrightarrow{\map{E}}(\Alg{B},\beta)$. First note that 
\[
\retro^{\Petz\;\;\;*}_{\alpha,\map{E}}(A)=(\Ad_{\hat{\beta}^{1/2}}\circ\map{E}\circ\Ad_{\alpha^{1/2}})(A)+\Tr[P_{\alpha}A]\frac{P_{\beta}^{\perp}}{\Tr[P_{\alpha}]}
\]
for all $A\in\Alg{A}$
by definition of the Hilbert--Schmidt adjoint. Hence, 
\[
\begin{split}
\big(\retro^{\Petz}_{\beta,\retro^{\Petz}_{\alpha,\map{E}}}\circ\Ad_{P_{\alpha}}\big)(A)
&=\big(\Ad_{\beta^{1/2}}\circ\retro^{\Petz\;\;\;*}_{\alpha,\map{E}}\circ\Ad_{\hat{\alpha}^{1/2}}\big)(A)
=(\Ad_{P_{\beta}}\circ\map{E}\circ\Ad_{P_{\alpha}})(A)\\
&=(\Ad_{P_{\alpha}}\circ\map{E}^*\circ\Ad_{P_{\beta}})^*(A)
=(\Ad_{P_{\alpha}}\circ\map{E}^*)^*(A)
=(\map{E}\circ\Ad_{P_{\alpha}})(A),
\end{split}
\]
where Lemma~\ref{lem:nullspaceE} was used in the second-last equality.
\eprf

\bq
Is $\retro^{\Petz}$ the unique retrodiction functor for $\NCFinStoch$? 
\eq

\bq
What categorical assumptions on a semicartesian category guarantee the existence of a retrodiction functor in such a way that both $\FinStoch$ and $\NCFinStoch$ satisfy such assumptions?
\eq

\bq
If a semicartesian category admits a retrodiction functor, what categorical assumptions guarantee that it is unique?
\eq

\bigskip
\noindent
{\bf Acknowledgements.}
I thank Francesco Buscemi, Tobias Fritz, Nicholas Gauguin Houghton-Larsen, Tom{\'a}{\v s} Gonda, and Paolo Perrone for discussions. 
I also thank Noson Yanofsky and Mahmoud Zeinalian for the invitations to speak about this topic at the New York City Category Theory Seminar and the Topology, Geometry, and Physics Seminar, respectively, at the CUNY Graduate Center, and for the encouraging discussions with the seminar participants.
This work was partially supported by MEXT-JSPS Grant-in-Aid for Transformative Research Areas (A) ``Extreme Universe'', No.\ 21H05183.
This work also acknowledges support from the Blaumann Foundation.

\addcontentsline{toc}{section}{\numberline{}Bibliography}
\bibliographystyle{eptcs}
\bibliography{references}

\Addresses

\end{document}